\documentclass[11pt]{article}
\usepackage{amssymb,amsmath,pstricks}
\makeatletter\@addtoreset {equation}{section}\makeatother

\newtheorem{theorem}{Theorem}
\newtheorem{lemma}{Lemma}
\newtheorem{remark}{Remark}
\newtheorem{corollary}{Corollary}
\newtheorem{definition}{Definition}

\newtheorem{example}{Example}
\newtheorem{proposition}{Proposition}

\newenvironment{proof}{
    \noindent {\it Proof.}}{\hfill$\Box$
}

\usepackage[dvips]{epsfig}
\usepackage{graphicx}

\begin{document}

\title{\bf Coupled-mode equations and gap solitons in a two-dimensional
nonlinear elliptic problem with a separable periodic potential}

\author{Tom\'{a}\v{s} Dohnal$^1$\footnote{Current address: Institut f\"{u}r Angewandte und Numerische Mathematik,
Universit\"{a}t Karlsruhe, Germany} , Dmitry Pelinovsky$^2$, and Guido Schneider$^3$ \\
{\small $^1$ Seminar for Applied Mathematics, ETH Z\"{u}rich,
Switzerland} \\ {\small $^2$ Department of Mathematics, McMaster University,
Hamilton, Ontario, Canada} \\
{\small $^3$ Institut f\"{u}r Analysis, Dynamik und Modellierung,
Universit\"{a}t Stuttgart, Germany} }

\date{\today}
\maketitle

\begin{abstract}
We address a two-dimensional nonlinear elliptic problem with a
{\em finite}-amplitude periodic potential. For a class of {\em
separable symmetric} potentials, we study the bifurcation of the
{\em first band gap} in the spectrum of the linear Schr\"{o}dinger
operator and the relevant {\em coupled-mode equations} to describe
this bifurcation. The coupled-mode equations are derived by the
rigorous analysis based on the Fourier--Bloch decomposition and
the Implicit Function Theorem in the space of bounded continuous
functions vanishing at infinity. Persistence of {\em reversible
localized} solutions, called {\em gap solitons}, beyond the
coupled-mode equations is proved under a non-degeneracy assumption
on the kernel of the linearization operator. Various branches of
reversible localized solutions are classified numerically in the
framework of the coupled-mode equations and convergence of the
approximation error is verified. Error estimates on the
time-dependent solutions of the Gross--Pitaevskii equation approximated by solutions
of the coupled-mode equations are obtained for a finite-time interval.
\end{abstract}

\section{Introduction}

Interplay between nonlinearity and periodicity is the focus of recent studies
in different branches of nonlinear physics and applied mathematics. Physical
applications of nonlinear systems with periodic potentials range from
nonlinear optics, in the dynamics of guided waves in inhomogeneous optical
structures and photonic crystal lattices, to atomic physics, in the dynamics
of Bose--Einstein condensate droplets in periodic potentials, and from
condensed matter, in Josephson-junction ladders, to biophysics, in various
models of the DNA double strand. The paramount significance for these models
is the possibility of spatial localization, that is emergence of nonlinear localized
structures residing in the spectral band gaps of the periodic potentials.
To describe this phenomenon, the primary equations of physics are typically simplified to
the Gross--Pitaevskii equation, which we shall study in our article in the space of two
dimensions. More precisely, we consider the two-dimensional Gross--Pitaevskii equation
in the form
\begin{equation}
\label{GP} i E_t = - \nabla^2 E + V(x) E + \sigma |E|^2 E,
\end{equation}
where $E(x,t) : \mathbb{R}^2 \times \mathbb{R} \mapsto
\mathbb{C}$, $\nabla^2 = \partial_{x_1}^2 + \partial_{x_2}^2$,
$V(x) : \mathbb{R}^2 \mapsto \mathbb{R}$, and $\sigma = \pm 1$.
Time-periodic solutions of the Gross--Pitaevskii equation are found
from the solutions of the nonlinear elliptic problem
\begin{equation}
\label{stationary} \nabla^2 \phi(x) + \omega \phi(x) = V(x) \phi +
\sigma |\phi(x)|^2 \phi(x),
\end{equation}
where $\phi(x) : \mathbb{R}^2 \mapsto \mathbb{C}$ and $\omega \in
\mathbb{R}$ arise in the substitution $E(x,t) = \phi(x) e^{-i
\omega t}$. It is known that localized solutions of the elliptic
problem (\ref{stationary}) with a periodic potential $V(x)$, called
{\em gap solitons}, exist in {\em every} finite gap of the spectrum
of the Schr\"{o}dinger operator $L = -\nabla^2 + V(x)$ and in the
semi-infinite gap for $\sigma = -1$ \cite{Stuart,Pankov}. Bifurcations of
localized solutions from edges of the spectral bands were
studied earlier in \cite{heinz,kupper}.

Coupled-mode equations were used by physicists for the analysis of
existence, stability and dynamics of gap solitons \cite{mills,SS}. A
justification of the one-dimensional coupled-mode equations in the
context of the elliptic problem (\ref{stationary}) with $x\in
\mathbb{R}$ was carried out in \cite{PSn07} in the limit of {\em
small-amplitude} periodic potentials $V(x)$. (A justification of
time-dependent coupled-mode equations on finite-time intervals
was done earlier in \cite{goodman,SU}.) In the limit of small-amplitude
potentials, narrow gaps of the spectrum of the Schr\"{o}dinger operator
$L = - \partial_x^2 + V(x)$ bifurcate from resonant points of the spectrum of $L_0 =
- \partial_x^2$, while the Bloch modes of $L$
bifurcate from the Fourier modes of $L_0$. In other words,
photonic band gaps open generally for small-amplitude
one-dimensional periodic potentials $V(x)$. {\em Small-amplitude}
gap solitons of the elliptic problem (\ref{stationary}) in $x \in
\mathbb{R}$ reside in the narrow gaps of $L$ for $V(x) \neq 0$
according to the approximation obtained from the coupled-mode equations \cite{PSn07}.

We refer to the opening of a spectral gap under a small change of the potential $V(x)$
as to the bifurcation of the band gap. Bifurcations of band gaps do not occur for
small-amplitude multi-dimensional periodic potentials. This is
caused by an overlap of the spectral bands of $L_0 = -\nabla^2$ in
the first Brillouin zone if $x \in \mathbb{R}^N$ and $N \geq 2$
\cite{Kuchment}. As a result, photonic band gaps in
multi-dimensional potentials open only at some finite amplitudes
of the periodic potential $V(x)$ and the resonant eigenfunctions
are given by the Bloch modes of $L = -\nabla^2 + V(x)$ rather than
by the Fourier modes of $L_0$. Although the coupled-mode equations
were also derived for multi-dimensional problems with small
periodic potentials \cite{ACD95,AFI04,AP05,AJ98,DA05} and the
resonant Fourier modes were used for the approximation of the full
solution, the applicability of these coupled-mode equations
remains an open issue for a rigorous analysis. Bloch mode
decomposition has been also used in one dimension for
finite-amplitude periodic potentials to derive coupled-mode
equations \cite{dSSS96}. The corresponding unperturbed one-dimensional
potential, however, has to be of a special type to admit a finite number of
open gaps, such that a new gap is opened under a small
perturbation.

In this paper we derive coupled-mode equations for wavepackets in
narrow band gaps of a {\em finite-amplitude} periodic potential by
using the Fourier--Bloch decomposition and the rigorous analysis
based on the Implicit Function Theorem in the space of bounded
continuous functions vanishing at infinity. The coupled-mode
equations we derive here take the form of coupled nonlinear
Schr\"{o}dinger (NLS) equations. These equations differ from the
first-order coupled-mode equations
exploited earlier \cite{SS}. Similar coupled NLS
equations have been recently
derived in \cite{SY} near band edges of the well-separated
spectral bands and in \cite{Konotop} for tunnelling problems.
Unlike these works relying on numerical approximations, we justify
the derivation of the coupled-mode equations and prove the
persistence of localized solutions in the full nonlinear problem
(\ref{stationary}). Although details of our analysis are given
only for the bifurcation of the first band gap in the spectrum of
$L$, a similar analysis can be developed for bifurcations of other
band gaps and for bifurcations of the localized solutions near
band edges of the well-separated spectral bands.

Our derivation is developed for the class of
{\em separable} potentials
\begin{equation}
\label{separable-potential} V(x_1,x_2) = \eta \left[ W(x_1) +
W(x_2) \right], \quad \eta \in \mathbb{R},
\end{equation}
where the function $W(x)$ is assumed to be real-valued, bounded,
piecewise-continuous, and $2\pi$-periodic on $x \in \mathbb{R}$.
To simplify the details of our analysis, we assume that $W(-x) = W(x)$ on $x
\in \mathbb{R}$, such that the solution set of the elliptic
problem (\ref{stationary}) includes functions satisfying
one of the following two \textit{reversibility} constraints
\begin{equation}\label{reversibility}
\phi(x_1,x_2) = s_1 \bar{\phi}(-x_1,x_2) = s_2 \bar{\phi}(x_1,-x_2)
\end{equation}
or
\begin{equation}
\label{reversibility-2} \phi(x_1,x_2) = s_1 \bar{\phi}(x_2,x_1) =
s_2 \bar{\phi}(-x_2,-x_1),
\end{equation}
where $s_1, s_2 = \pm 1$.

Our strategy is to show that there may exist a value $\eta =
\eta_0$, for which the first band gap opens due to the resonance
of three lowest-order Bloch modes for three spectral bands of the
operator $L = -\nabla^2 + V(x)$. The new small parameter $\epsilon
:= \eta - \eta_0$ is then used for the bifurcation theory of Bloch
modes which results in the {\em algebraic} coupled-mode equations
for nonlinear interaction of the three resonant modes. The
Fourier--Bloch decomposition is used for the approximation of
localized solutions and for the derivation of the {\em
differential} coupled-mode equations with the second-order
derivative terms. The main idea behind our technique is that
a differential operator after the Fourier--Bloch decomposition
becomes a pseudo--differential operator whose symbol is the multiple-valued
dispersion relation. Then only the relevant (three) branches
of the resonant modes can be taken into account, which leads to a coupled-mode system.

If the linearization operator of the coupled NLS
equations at the localized solutions is non-degenerate, the
persistence of the localized symmetric solutions in the full
nonlinear problem (\ref{stationary}) is proved with the Implicit
Function Theorem. Localized symmetric solutions of the coupled NLS
equations are approximated numerically and the convergence rate
for the error of approximation is studied. Finally, we study
time-dependent localized solutions of the Gross--Pitaevskii
equation and the coupled-mode system and control smallness of the
distance between the two solutions on a finite-time interval.

The article is structured as follows. Section \ref{S:Strum_Liou}
reviews elements of the Sturm--Liouville theory for the separable
potentials. Section \ref{S:algeb_CME} contains a derivation of the
algebraic coupled-mode equations for three resonant Bloch modes.
Section \ref{S:diff_CME} gives details of the projection technique
for the derivation of the differential coupled-mode equations for
localized solutions. The persistence of localized reversible
solutions under a non-degeneracy assumption on the linearization
operator is proved in Section \ref{S:reversible_GS}. Numerical
approximations of the localized solutions, the associated
linearization operators, and the convergence of the approximation
error are obtained in Section \ref{S:numerical}. Section
\ref{S:t_dep_CME} extends the results to the time-dependent case
for finite time intervals. Section \ref{S:generalizations}
discusses relevant generalizations.

\section{Sturm--Liouville theory for separable potentials}\label{S:Strum_Liou}

\begin{figure}[htpb]
\centering \setlength{\unitlength}{1cm}
\begin{picture}(5.5,4.5)(-1.5,-0.5)
  \put(1.6,1.0){\huge $B_0$}
      \put(0.65,3.2){$X'$}
  \pspolygon[linecolor=black,linewidth=0.04cm,fillstyle=vlines,fillcolor=lightgray](0.5,0.5)(3,0.5)(3,3)(0.5,0.5)
  \psline(-.3,.5)(3.9,.5)  
   \psline(0.5,-.3)(0.5,3.7)  
  \psline{->}(3.9,.5)(4.1,.5)
  \psline{->}(.5,3.7)(.5,3.9)
  \put(-0.1,3.7){$k_2$}
   \put(4.1,0.1){$k_1$}
   \psline(0.5,3)(3,3)
   \put(0.1,0){$\Gamma$}
   \put(3.2,0.7){$X$}
   \put(3.2,3.2){$M$}
      \put(0.1,2.85){$\frac{1}{2}$}
   \put(2.85,0.08){$\frac{1}{2}$}
\end{picture}
\caption{The first irreducible Brillouin zone $B_0$ for a
two-dimensional separable potential.} \label{F:B_ir}
\end{figure}
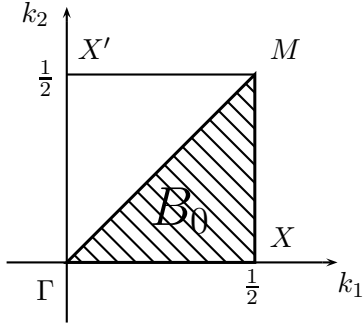

It is typically expected that the band gaps in the spectrum of the linear Schr\"{o}dinger
operator with a two-dimensional periodic potential open
at the extremal values of the Bloch band surfaces \cite{Kuchment}.
For a general potential, however, the extremal values may
occur anywhere within the first irreducible Brillouin zone $B_0$
in the quasi-momentum space $(k_1,k_2)$ \cite{HKSW_prep}. We will show here
that the extremal values for a separable
potential occur only at the vertex points $\Gamma$, $X$ and $M$ on
the boundary $\partial B_0$. Figure \ref{F:B_ir} shows the
irreducible Brillouin zone $B_0$ and the vertex points for a two-dimensional
separable potential $V(x)$.

Let us consider the spectral problem for the linear
Schr\"{o}dinger operator associated with the separable periodic
potential (\ref{separable-potential}):
\begin{equation}
\label{lin_eval} -\nabla^2 u(x_1,x_2) + \eta \left[ W(x_1) +
W(x_2) \right] u(x_1,x_2) = \omega u(x_1,x_2).
\end{equation}
By using the separation of variables $u(x_1,x_2) = f_1(x_1) f_2(x_2)$
and the parametrization $\omega = \xi_1 + \xi_2$, we obtain the
uncoupled eigenvalue problems
\begin{equation}
\label{separated-ODE}
-f_j''(x_j) + \eta W(x_j) f_j(x_j) = \xi_j f_j(x_j), \ j=1,2.
\end{equation}
The Bloch modes and the band surfaces for the one-dimensional problems
(\ref{separated-ODE}) are introduced according to the regular
Sturm--Liouville problem
\begin{equation}\label{E:Bloch_eval_problem}
\left\{ \begin{array}{l} -u''(x) + \eta W(x) u(x) = \rho u(x),
\;\; 0 \leq x \leq 2\pi, \\ u(2\pi) = e^{i 2\pi k}u(0),
\end{array} \right.
\end{equation}
where $k$ is a quasi-momentum defined on the interval
$\mathbb{T} = \left[ -\frac{1}{2},\frac{1}{2}\right ]$. The eigenfunctions of
the Sturm--Liouville problem (\ref{E:Bloch_eval_problem}) are periodic
with respect to $k$ with the period one. By Theorem 2.4.3 in
\cite{Eastham}, there exists a countable infinite set of
eigenvalues $\{ \rho_n(k) \}_{n \in \mathbb{N}}$ for each $k \in
\mathbb{T}$, which can be ordered as
$$
\rho_1(k) \leq \rho_2(k) \leq \rho_3(k) \leq ...
$$
By Theorem 4.2.3 in \cite{Eastham}, if $W(x)$ is a bounded and
$2\pi$-periodic potential, eigenvalues
$\{ \rho_n(k) \}_{n \in \mathbb{N}}$ have a uniform asymptotic distribution on $k \in
\mathbb{T}$, such that
\begin{equation}
\label{asymptotic-distribution} C_- n^2 \leq |\rho_n(k)| \leq C_+
n^2, \qquad \forall n \in \mathbb{N}, \;\; \forall k \in
\mathbb{T},
\end{equation}
for some constants $C_{\pm} > 0$.

Let $\{ u_n(x;k) \}_{l \in \mathbb{N}}$ be the corresponding set
of eigenfunctions of the Sturm--Liouville problem
(\ref{E:Bloch_eval_problem}), such that $u_n(x+2\pi;k) = u_n(x;k)
e^{2 \pi i k}$. By Theorem XIII.89 in \cite{RS}, if $W(x)$ is a
bounded, piecewise-continuous and $2\pi$-periodic potential, then
the eigenvalue $\rho_n(k)$ and the Bloch function $u_n(x;k)$
are analytic in $k \in \mathbb{T} \backslash \left\{ 0,
\pm\frac{1}{2} \right\}$ and continuous at the points $k = 0$
and $k = \pm \frac{1}{2}$. By Theorem XIII.95 in \cite{RS}, the
eigenvalue $\rho_n(k)$ is extended on a smooth Riemann surface
in the neighborhood of the points $k = 0$ and $k = \frac{1}{2}$ if
the $n^{\rm th}$ spectral band is disjoint from the adjacent
$(n\pm 1)^{\rm th}$ spectral bands by a non-empty band gap.

By Theorem XIII.90 in \cite{RS}, if $W(x)$ is a bounded,
piecewise-continuous and $2\pi$-periodic potential, then the
spectrum of $L_{1D} = -
\partial_x^2 + \eta W(x)$ in $L^2(\mathbb{R})$ is
absolutely continuous and consists of the union of the intervals
in the range of the functions $\rho_n(k)$ on $k \in \mathbb{T}$ for $n
\in \mathbb{N}$, where $\rho_n(-k) = \rho_n(k)$. Moreover, the
extremal values of $\rho_n(k)$ occur only at the points $k =
0$ and $k = \pm \frac{1}{2}$. These points correspond to an alternating
sequence of the maximum and minimum values of the eigenvalues $\{ \rho_n(k)
\}_{n \in \mathbb{N}}$ according to the following formula
$$
\arg \min_{k \in \mathbb{S}} \rho_{2m-1}(k) = 0, \;  \arg
\max_{k \in \mathbb{S}} \rho_{2m-1}(k) = \pm \frac{1}{2}, \;
\arg \min_{k \in \mathbb{S}} \rho_{2m}(k) = \pm \frac{1}{2}, \;
\arg \max_{k \in \mathbb{S}} \rho_{2m}(k) = 0,
$$
for all $m \in \mathbb{N}$. The Bloch function
$u_n(x;k)$ is $2\pi$-periodic if $k = 0$ and $2 \pi$-antiperiodic
if $k = \pm \frac{1}{2}$.  Let us denote the periodic and
antiperiodic eigenfunctions and the corresponding eigenvalues by
\begin{equation}
\label{eigenvalue-periodic} \psi_n(x) = u_n(x;0), \;\; \lambda_n =
\rho_n(0) \quad \text{and} \quad \varphi_n(x) =
u_n\left(x;\pm\frac{1}{2}\right), \;\; \mu_n = \rho_n\left(\pm
\frac{1}{2}\right), \quad n \in \mathbb{N}.
\end{equation}
By Theorems 2.3.1 and 3.1.2 in \cite{Eastham}, these eigenvalues are ordered by
$\lambda_1 < \lambda_2 \leq \lambda_3 \leq ...$ and $\mu_1 < \mu_2
\leq \mu_3 \leq ...$, while the corresponding eigenfunctions
$\psi_n(x)$ and $\varphi_n(x)$ have precisely $n-1$ zeros (nodes)
on the interval $(-\pi,\pi)$. If $W(-x) = W(x)$, the
eigenfunctions $\psi_n(x)$ and $\varphi_n(x)$ are even for odd $n$
and odd for even $n$. Each set of the eigenfunctions $\{ \psi_n(x) \}_{n \in
\mathbb{Z}}$ and $\{ \varphi_n(x) \}_{n \in \mathbb{Z}}$ is orthogonal in
$L^2([-\pi,\pi])$.

Eigenvalues of the same one-dimensional operator $L_{1D} = - \partial_x^2 + \eta
W(x)$ for $4 \pi$-periodic eigenfunctions consist of the union of
the eigenvalues $\{ \lambda_n \}_{n \in \mathbb{N}}$ and $\{ \mu_n
\}_{n \in \mathbb{N}}$. Since there are at most two eigenvalues of
the second-order operator $L_{1D}$ and by Theorem 2.3.1 of
\cite{Eastham}, we obtain that
\begin{equation}
\label{ordering} \lambda_1 < \mu_1 < \mu_2 < \lambda_2 \leq
\lambda_3 < \mu_3 \leq \mu_4 < \lambda_4 \leq \lambda_5 < ...
\end{equation}
In particular, the first band gap is always non-empty for a
non-constant potential $W(x)$ (see Theorem XIII.91(a) in
\cite{RS}). Due to this ordering, the lowest value of $\omega = \rho_{n_1}(k_1) +
\rho_{n_2}(k_2)$ for the crossing of the Bloch band surfaces associated
with the two-dimensional separable potential $V(x_1,x_2)$ occurs at $\omega = \omega_0
= \lambda_1 + \mu_2 = 2 \mu_1$. Using these facts, we
obtain the following results.

\begin{lemma}
Extremal values of the Bloch band surfaces $\omega =
\rho_{n_1}(k_1) + \rho_{n_2}(k_2)$ for $(k_1,k_2)\in B_0$ occur only at the vertex
points of the boundary $\partial B_0$.
\end{lemma}

\begin{proof}
Because $\omega = \rho_{n_1}(k_1) + \rho_{n_2}(k_2)$ for a
separable potential (\ref{separable-potential}), we have $\nabla \omega =
\left[ \rho_{n_1}'(k_1),\rho_{n_2}'(k_2) \right]^T$, where
$\rho_n(k)$ are analytic on $k \in \mathbb{T}\backslash
\left\{ 0, \pm \frac{1}{2} \right\}$. If an extremal point occurs
in the interior of the Brillouin zone $B_0$, then $\nabla \omega =
0$ at $k = k_0 \in B_0$. However, $\rho_n'(k) \neq 0$ for any $0 <
|k| < \frac{1}{2}$ and any $n \in \mathbb{N}$ by Theorem XIII.90
in \cite{RS}. Similarly, the extremal points cannot occur in the
interior of the boundary $\partial B_0$. Therefore, the extremal
values of $\omega$ occur only at the vertex points $\Gamma$, $X$
and $M$ on $\partial B_0$.
\end{proof}

\begin{remark}
{\rm The derivative $\rho_n'(k)$ may be non-zero at $k = 0$ and $k
= \frac{1}{2}$ if the $n^{\rm th}$ spectral band touches the
adjacent $(n\pm 1)^{\rm th}$ spectral bands. This happens
when the corresponding eigenvalue $\lambda_n$ or $\mu_n$ is double
degenerate with the equality sign in the ordering
(\ref{ordering}). However, $\lambda_1$, $\mu_1$ and $\mu_2$ are
always simple and, therefore, $\nabla \omega = 0$ at least for the
first three spectral bands $\omega = \rho_{n_1}(k_1) +
\rho_{n_2}(k_2)$ at $\omega = \omega_0 = \lambda_1 +
\mu_2 = 2 \mu_1$.}
\end{remark}

\begin{lemma}
Assume that the resonant condition $\lambda_1 + \mu_2 = 2 \mu_1
\equiv \omega_0$ for the $2\pi$-periodic eigenvalue $\lambda_1$ and the $2\pi$-antiperiodic
eigenvalues $\mu_1, \mu_2$ of $L_{1D} = -\partial_x^2 + \eta
W(x)$ is satisfied for $\eta = \eta_0$. Then there are exactly three
resonant Bloch modes at $\omega = \omega_0$ and $\eta = \eta_0$ in
the spectral problem (\ref{lin_eval}):
$$
\Phi_1 = \psi_1(x_1) \varphi_2(x_2), \;\; \Phi_2 = \varphi_2(x_1)
\psi_1(x_2), \;\; \Phi_3 = \varphi_1(x_1) \varphi_1(x_2).
$$
The three resonant modes are orthogonal to each other with respect
to the $4\pi$-periodic inner product
$$
(f,g) = \int_{-2\pi}^{2\pi}\int_{-2\pi}^{2\pi} \bar{f}(x_1,x_2) g(x_1,x_2) dx_1 dx_2.
$$
\label{lemma-resonance}
\end{lemma}

\begin{proof}
There are exactly three resonant modes for {\em any} separable
potential (\ref{separable-potential}) because $\lambda_1$ and
$\mu_1$ are the smallest non-degenerate eigenvalues of the operator
$L_{1D}$ and the double degeneracy $\lambda_1 + \mu_2 = \mu_2 +
\lambda_1$ is due to the symmetry with respect to the interchange of
the variables $x_1$ and $x_2$. The eigenfunctions $\{ \psi_n(x) \}_{n
\in \mathbb{N}}$ and $\{ \varphi_n(x) \}_{n\in\mathbb{N}}$ are all
orthogonal to each other in $L^2([-2\pi,2\pi])$, which leads to
the orthogonality of the three resonant modes.
\end{proof}

\begin{remark}
{\rm We note that the eigenfunctions $\psi_1(x)$ and
$\varphi_1(x)$ are not orthogonal on the interval $[-\pi,\pi]$
since both of them are positive by Theorem 3.1.2 in
\cite{Eastham}. However, a linear combination of $\psi_1(x)$ and
$\varphi_1(x)$ belongs to the class of $4 \pi$-periodic functions
and the two eigenfunctions are orthogonal on the double-length
interval $[-2\pi,2\pi]$.}
\end{remark}

\begin{remark}
{\rm By Theorem 6.10.5 in \cite{Eastham}, there are finitely many
gaps in the spectral problem (\ref{lin_eval}) with a separable
potential. (In fact, there are finitely many gaps for general smooth periodic potentials in dimension 2 or higher \cite{Parnov_08}.)
Lemma \ref{lemma-resonance} indicates a
bifurcation when one of these gaps at the lowest value of $\omega$
may open. It is observed for many examples of the separable
potential (\ref{separable-potential}) that this bifurcation leads
to the {\em first} band gap in the spectrum of $L = - \nabla^2 +
\eta [W(x_1) + W(x_2)]$ if $\eta$ is increased from $\eta = 0$. }
\end{remark}

\begin{figure}[htbp]
\begin{center}
\includegraphics[height=6cm]{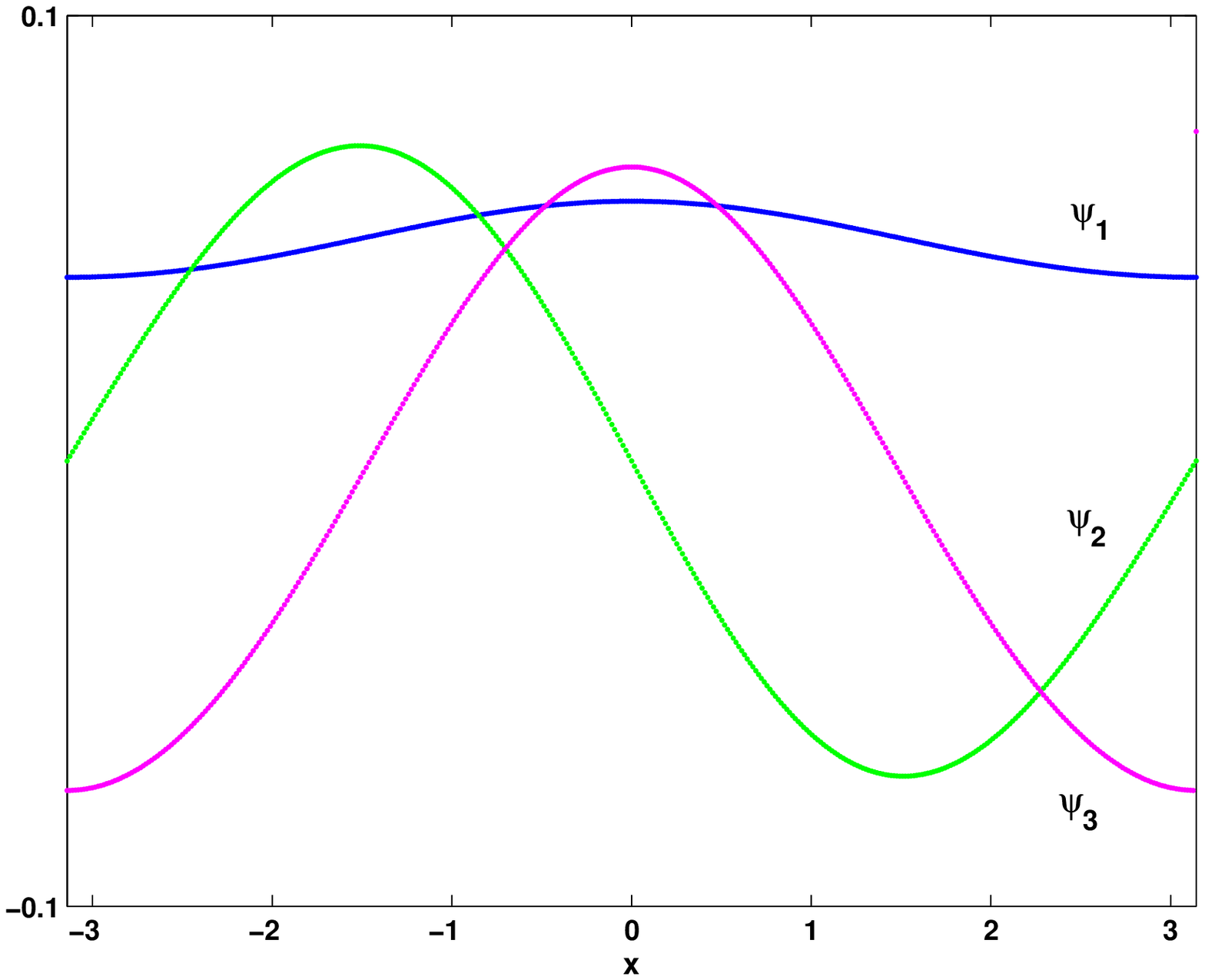}
\includegraphics[height=6cm]{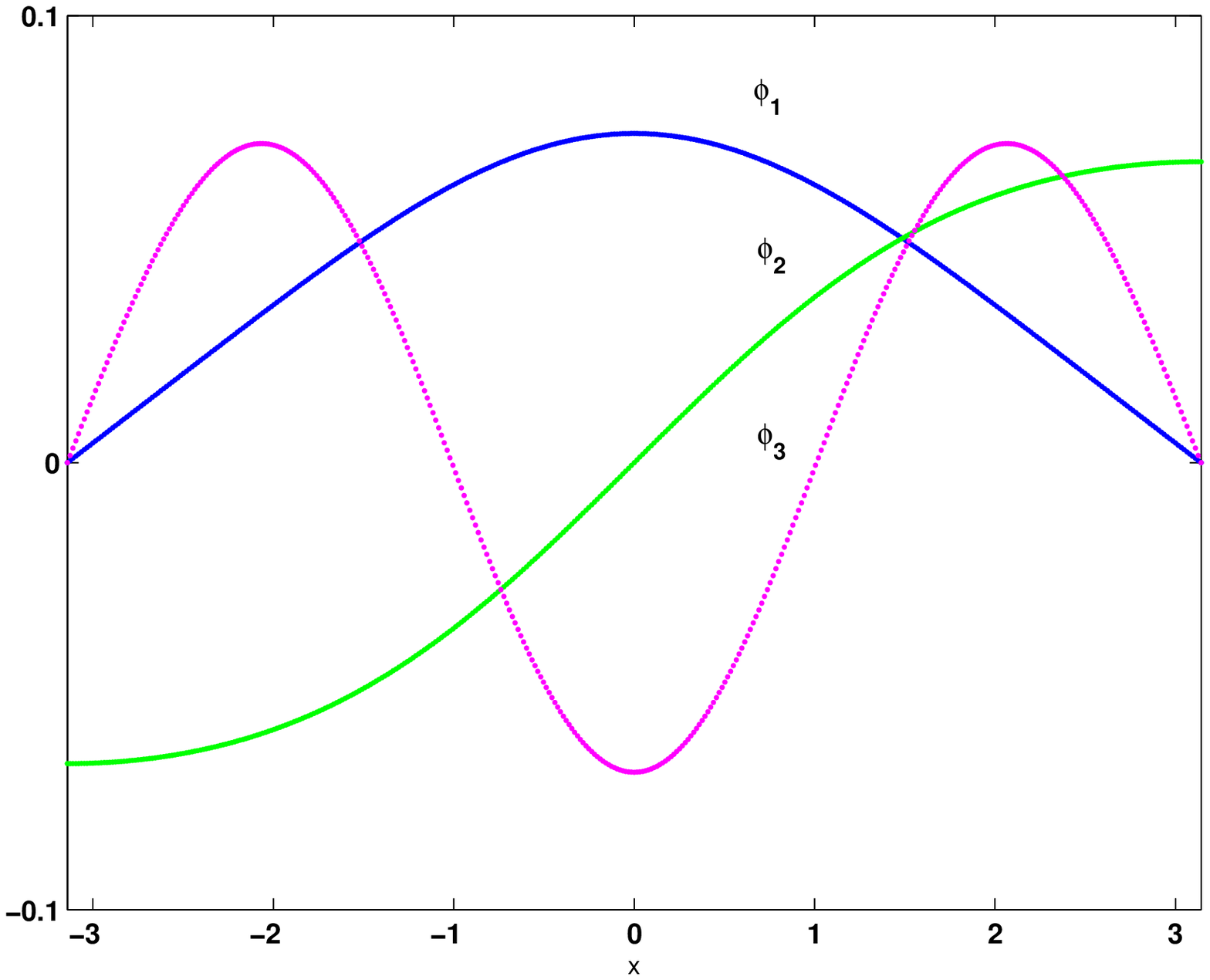}
\end{center}
\caption{The $2\pi$-periodic (left) and $2\pi$-antiperiodic
(right) eigenfunctions of $L_{1D} = -\partial_x^2 + \eta_0 W(x)$
with $W(x) = 1 - \cos x$.} \label{F:eigenfunctions}
\end{figure}

\begin{figure}[htbp]
\begin{center}
\includegraphics[height=6cm]{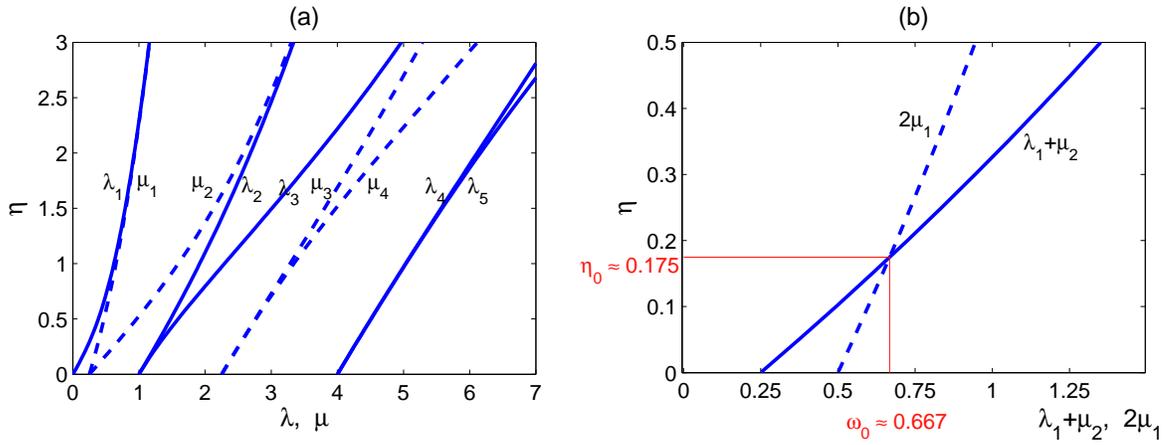}
\end{center}
\caption{(a) Eigenvalues of $L_{1D} = -\partial_x^2 + \eta (1 -
\cos x)$ for $2\pi$-periodic (solid) and $2\pi$-antiperiodic
(dashed) boundary conditions. (b) The lowest resonance occurs at
$\omega = \omega_0$ and $\eta = \eta_0$ when $\lambda_1 +\mu_2 =
2\mu_1$.} \label{F:bands}
\end{figure}

\begin{figure}[htbp]
\begin{center}
\includegraphics[height=6cm]{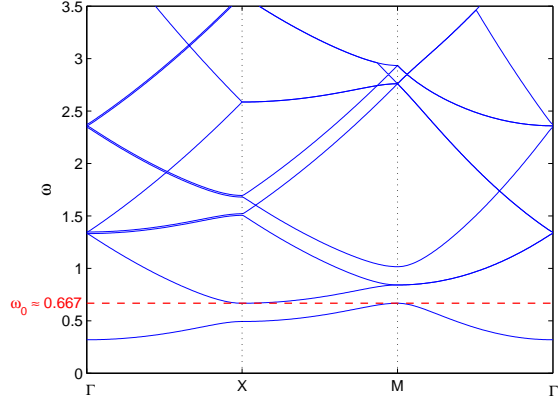}
\end{center}
\caption{Band diagram along $\partial B_0$  for the spectral
problem (\ref{lin_eval}) with $W(x) = 1 - \cos x$ and
$\eta=\eta_0$.} \label{F:2D_band_diagram}
\end{figure}

\begin{example}
\label{example-1} {\rm Let $W(x) = 1 - \cos x$. Numerical
approximations of the eigenvalues of the Sturm--Liouville problem
(\ref{E:Bloch_eval_problem}) are computed with the use of the
second-order central difference method. The eigenfunctions $\psi_n(x)$
and $\varphi_n(x)$ are plotted in Figure \ref{F:eigenfunctions}
for $n = 1,2,3$, while the dependence of the first eigenvalues
$\lambda_n$ and $\mu_n$ on $\eta$ is plotted in Figure
\ref{F:bands}(a). Figure \ref{F:bands}(b) shows the intersection
of $2\mu_1$ and $\lambda_1 + \mu_2$ at the bifurcation value $\eta
= \eta_0 \approx 0.1745$, when $\lambda_1 \approx 0.1595$, $\mu_1
\approx 0.3336$ and $\mu_2 \approx 0.5077$, such that $\omega_0 =
\lambda_1 + \mu_2 = 2 \mu_1 \approx 0.6672$. The perturbation
behavior shown in Figure \ref{F:bands} was studied analytically in \cite{arnold}.
Figure
\ref{F:2D_band_diagram} illustrates the band structure along
$\partial B_0$ of the full spectral problem (\ref{lin_eval}) with
$\eta = 0.1745$ clearly revealing the resonance at $\omega \approx
0.6672$ and the two resonant modes at the points $X$ and
$M$. The third resonant mode lies at $X'$ because the symmetry $V(x_1,x_2)=V(x_2,x_1)$ implies
$\omega(k_1,k_2)=\omega(k_2,k_1)$.}
\end{example}

In what follows, we consider the bifurcation of non-trivial
solutions of the nonlinear elliptic problem (\ref{stationary})
with the separable potential (\ref{separable-potential}) in the
lowest band gap described by Lemma \ref{lemma-resonance}. Let
$\epsilon = \eta - \eta_0$, $\omega = \omega_0 + \epsilon \Omega$,
$\phi(x) = \sqrt{\epsilon} \Phi(x)$ and rewrite the nonlinear
elliptic problem (\ref{stationary}) in the form
\begin{equation}
\label{bifurcation} {\cal L}_0 \Phi(x) = \epsilon \Omega \Phi(x) -
\epsilon \left[ W(x_1) + W(x_2) \right] \Phi(x) - \epsilon
\sigma |\Phi(x)|^2 \Phi(x),
\end{equation}
where $\Phi : \mathbb{R}^2 \mapsto \mathbb{C}$ and ${\cal L}_0 = -
\nabla^2 + \eta_0 \left[ W(x_1) + W(x_2) \right] - \omega_0$. Two
classes of non-trivial solutions of the bifurcation problem
(\ref{bifurcation}) are considered for small $\epsilon$: bounded
$4\pi$-periodic solutions (Section \ref{S:algeb_CME}) and bounded
decaying solutions (Sections \ref{S:diff_CME}--\ref{S:numerical}).

\section{Algebraic coupled-mode equations}\label{S:algeb_CME}

We consider here bounded periodic solutions $\Phi(x)$ of the
bifurcation problem (\ref{bifurcation}) for small $\epsilon$. Our
results depend on the period of the periodic solutions $\Phi(x)$.
Since the potential $V(x)$ is separable, the $2\pi$-periodic
($2\pi$-antiperiodic) function $\Phi(x)$ can be represented by the
series of Bloch modes associated with the eigenvalue problem
(\ref{E:Bloch_eval_problem}) on $x \in \mathbb{R}$ for $k = 0$ ($k
= \pm \frac{1}{2}$), while the $4\pi$-periodic functions $\Phi(x)$
can be represented by the series of both $2\pi$-periodic and
$2\pi$-antiperiodic Bloch modes. In what follows, we shall use
symbols $L^2(\mathbb{P}_k)$, $H^s(\mathbb{P}_k)$ and $C_b^0(\mathbb{P}_k)$
to denote the corresponding spaces of functions $u(x)$ on the interval $[0,2\pi]$
satisfying the boundary conditions $u(2\pi) = e^{i 2 \pi k} u(0)$. We shall also
use symbol $\vec{\mbox{\boldmath $\phi$}}$ to denote the vector of elements of
the sequence $\{ \phi_n \}_{n \in \mathbb{N}}$.

\begin{proposition}
Let $W(x)$ be a bounded and $2\pi$-periodic function. Let $\{
\rho_n(k) \}_{n \in \mathbb{N}}$ and $\{ u_n(x;k) \}_{n \in
\mathbb{N}}$ be sets of eigenvalues and eigenfunctions of the
Sturm--Liouville problem (\ref{E:Bloch_eval_problem}) that depends
on $k \in \mathbb{T}$, such that
\begin{equation}
\label{orthogonality-1} \langle u_n(\cdot;k), u_{n'}(\cdot;k)
\rangle_{2\pi} = \delta_{n,n'}, \qquad \forall n,n' \in
\mathbb{N}, \;\; \forall k \in \mathbb{T},
\end{equation}
where $\langle f,g \rangle_{2\pi} = \int_0^{2 \pi} \bar{f}(x) g(x)
dx$ and $\delta_{n,n'}$ is the Kronecker symbol. For any fixed $k
\in \mathbb{T}$, the set of eigenfunctions is complete in
$L^2(\mathbb{P}_k)$, such that there exists a unique set of
coefficients $\{ \phi_n \}_{n \in \mathbb{N}}$ in the
decomposition
\begin{equation}
\label{decomposition-1} \forall \phi(x) \in L^2(\mathbb{P}_k) :
\qquad \phi(x) = \sum_{n \in \mathbb{N}} \phi_n u_n(x;k),
\end{equation}
given by $\phi_n = \langle u_n(\cdot;k),\phi \rangle_{2\pi}$.
\label{proposition-1}
\end{proposition}

\begin{proof}
The statement of the proposition follows by Theorem XIII.88 of
\cite{RS}.
\end{proof}

\begin{lemma}
\label{corollary-decomposition-1} Let $\phi(x)$ be defined by the
decomposition (\ref{decomposition-1}), where $\vec{\mbox{\boldmath
$\phi$}}$ belongs to the vector space $l^1_s(\mathbb{N})$ with the
norm
\begin{equation}
\label{vector-space-1} \| \vec{\mbox{\boldmath $\phi$}}
\|_{l^1_s(\mathbb{N})} = \sum_{n \in \mathbb{N}} (1+n)^s |\phi_n|
< \infty.
\end{equation}
If $s > \frac{1}{2}$, then $\phi(x)$ is a continuous function on
$x \in \mathbb{R}$ and $\phi(x) = \psi(x) e^{i k x}$ with
$\psi(x+2\pi) = \psi(x)$.
\end{lemma}

\begin{proof}
By the triangle inequality, we obtain
$$
\| \phi \|_{C^0_b(\mathbb{R})} \leq \sum_{n \in \mathbb{N}}
|\phi_n| \| u_n(\cdot;k) \|_{C^0_b(\mathbb{P}_k)}.
$$
By Sobolev's Embedding Theorem, there exists a $C > 0$ such that
$\| u_n(\cdot;k) \|_{C^0_b(\mathbb{P}_k)} \leq C \| u_n(\cdot;k)
\|_{H^s(\mathbb{P}_k)}$ for any $s > \frac{1}{2}$, $n \in
\mathbb{N}$ and $k \in \mathbb{T}$. Since $W(x)$ is a bounded
potential, the squared norm $\| u_n(\cdot;k)
\|^2_{H^s(\mathbb{P}_k)}$ is equivalent to the integral
$$
\int_{0}^{2\pi} \left| (c_{\eta} + L_{1D})^{s/2} u_n(x;k) \right|^2
dx = (c_{\eta} + \rho_n(k))^{s},
$$
where $c_{\eta} > -\eta \min_{x \in \mathbb{R}} W(x)$ and
$(c_{\eta} + L_{1D})^{s/2}$ is defined using the spectral family
associated with the Sturm--Liouville operator $L_{1D} =
-\partial_x^2 + \eta W(x)$ in $L^2(\mathbb{P}_k)$. By the
asymptotic distribution of eigenvalues
(\ref{asymptotic-distribution}), we obtain that
$$
\| \phi \|_{C^0_b(\mathbb{P}_k)} \leq C \sum_{n \in \mathbb{N}}
|\phi_n| \| u_n(\cdot;k) \|_{H^s(\mathbb{P}_k)} \leq \tilde{C}
\sum_{n \in \mathbb{N}} (1 + n)^s |\phi_n|
$$
for some $\tilde{C} > 0$ and any $s > \frac{1}{2}$.
\end{proof}

If the potential $\eta W(x)$ is continued with respect to the
parameter $\eta$, the perturbation theory for an eigenvalue
$\rho_n(k)$ depends on whether $\rho_n(k)$ is simple or multiple. The following lemma
is a trivial statement of the perturbation theory, so we omit its proof.

\begin{lemma}
Let $W(x)$ be a bounded and $2\pi$-periodic function of $x \in
\mathbb{R}$. Fix $k \in \mathbb{T}$, $\eta_0 \in \mathbb{R}$, and
$n \in \mathbb{N}$. If $\rho_{n-1}(k) < \rho_n(k) <
\rho_{n+1}(k)$, then $\rho_n(k)$ and $u_n(x;k)$ depend
analytically on $\eta$ near $\eta = \eta_0$, such that
\begin{equation}
\label{leading-order}
\partial_{\eta} \rho_n(k) |_{\eta = \eta_0} =
\langle u_n(\cdot;k) |_{\eta = \eta_0}, W u_n(\cdot;k) |_{\eta =
\eta_0} \rangle_{2\pi}
\end{equation}
and
\begin{equation}
|\rho_n(k) - \rho_n(k) |_{\eta = \eta_0} -
\partial_{\eta} \rho_n(k) |_{\eta = \eta_0} (\eta - \eta_0) |
\leq C (\eta - \eta_0)^2  \label{Taylor-bound}
\end{equation}
for some $C > 0$ and sufficiently small $|\eta - \eta_0|$.
\label{lemma-perturbation-theory}
\end{lemma}

Due to the construction of the spectrum of the two-dimensional
spectral problem (\ref{lin_eval}), the results of Lemma
\ref{lemma-perturbation-theory} settle the question of the
splitting of the three resonant bands near the point $\omega =
\omega_0$ described in Lemma \ref{lemma-resonance}. Indeed, each
Bloch band surface (the graph of $\omega = \rho_{n_1}(k_1) + \rho_{n_2}(k_2)$)
corresponding to the three resonant bands includes a different
vertex point ($X'$, $X$, and $M$) at $\omega = \omega_0$ for
different values of $k$
in the set $\left\{\left(
0,\frac{1}{2} \right);  \left( \frac{1}{2},0 \right);  \left( \frac{1}{2},\frac{1}{2} \right)
\right\}$. Each eigenvalue of the spectral problem (\ref{lin_eval})
is simple and isolated if $k$ is fixed,
such that the continuation of each Bloch band surface in $\eta$ near
$\eta = \eta_0$ is determined by the expression
(\ref{leading-order}) and the separation of variables in $\omega =
\rho_{n_1}(k_1) + \rho_{n_2}(k_2)$. In other words, bifurcation of
the lowest band gap that occurs at $\omega = \omega_0$ can not be studied
if the function space is restricted to $L^2(\mathbb{P}_{k_1} \times \mathbb{P}_{k_2})$
for a fixed value of $(k_1,k_2)\in \mathbb{T}^2$. To detect the bifurcation, we shall
work in the function space $L^2_{\rm per}([-2\pi,2\pi]\times [-2\pi,2\pi])$,
where all three vertex points $X'$, $X$ and
$M$ correspond to the same $4\pi$-periodic boundary conditions on the function
$\Phi(x)$. According to Lemma \ref{lemma-resonance}, the three
resonant Bloch modes are orthogonal to each other in this space.

To work with the $4\pi$-periodic eigenfunctions, we shall
introduce a number of notations. Let $\{ \nu_n \}_{n \in
\mathbb{N}}$ and $\{ v_n(x) \}_{n \in \mathbb{N}}$ denote the sets
of eigenvalues and the corresponding $4 \pi$-periodic
eigenfunctions of the one-dimensional spectral problem
(\ref{E:Bloch_eval_problem}). These sets consist of the union of
eigenvalues and eigenfunctions defined by
(\ref{eigenvalue-periodic}). Suppose that the eigenvalues are
sorted in non-decreasing order, such that $\nu_1 = \lambda_1$,
$\nu_2 = \mu_1$, $\nu_3 = \mu_2$, $\nu_4 = \lambda_2$, $\nu_5 =
\lambda_3$ and so on, and that the eigenfunctions $v_1 = \psi_1,
v_2=\varphi_1, v_3=\varphi_2, v_4 = \psi_2, v_5= \psi_3$, etc.,
are normalized such that $\langle v_i,v_j\rangle_{4\pi} =
\delta_{i,j}$, where $\langle f,g\rangle_{4\pi} =
\int_{-2\pi}^{2\pi}\bar{f}(x)g(x)dx$. The statements of
Proposition \ref{proposition-1} and Lemma
\ref{corollary-decomposition-1} extend naturally to the new sets
of eigenfunctions in space $L^2_{\rm per}([-2\pi,2\pi])$.
To develop the nonlinear analysis of bifurcations of a non-trivial $4\pi$-periodic solution
$\Phi(x)$ on $x \in \mathbb{R}^2$, we shall study first the nonlinear vector
field $|\Phi(x)|^2 \Phi(x)$ acting on the
decomposition $\Phi(x) = \sum_{n \in \mathbb{N}} \phi_n v_n(x)$, where
$\vec{\mbox{\boldmath $\phi$}}$ belongs to the vector space $l^1_s(\mathbb{N})$.

\begin{lemma}
\label{lemma-vector-field} Let $W(x)$ be a bounded and
$2\pi$-periodic function. Let $\phi(x) = \sum_{n \in \mathbb{N}}
\phi_n v_n(x)$ and $|\phi(x)|^2 \phi(x) = \sum_{n \in \mathbb{N}}
g_n v_n(x)$.  If $0 < s < 1$, then there exists a $C > 0$ such that $\|
\vec{{\bf g}} \|_{l^1_s(\mathbb{N})} \leq C \| \vec{\mbox{\boldmath $\phi$}}
\|^3_{l^1_s(\mathbb{N})}$.
\end{lemma}

\begin{proof}
We only need to prove that the space $l^1_s(\mathbb{N})$ with $0 < s <
1$ forms a Banach algebra in the sense
\begin{equation}
\label{Banach-algebra} \forall \vec{\mbox{\boldmath $\phi$}},
\vec{\mbox{\boldmath $\varphi$}} \in l^1_s(\mathbb{N}) : \quad \|
\vec{\mbox{\boldmath $\phi$}} \star \vec{\mbox{\boldmath $\varphi$}}
\|_{l^1_s(\mathbb{N})} \leq C \| \vec{\mbox{\boldmath
$\phi$}}\|_{l^1_s(\mathbb{N}) } \|\vec{\mbox{\boldmath
$\varphi$}}\|_{l^1_s(\mathbb{N})},
\end{equation}
for some $C > 0$, where
$$
\left(\vec{\mbox{\boldmath $\phi$}} \star \vec{\mbox{\boldmath
$\varphi$}}\right)_n = \sum_{i \in \mathbb{N}} \sum_{j \in
\mathbb{N}} K_{n,i,j} \phi_{i} \varphi_{j}, \qquad \forall
n \in \mathbb{N}
$$
and
$$
K_{n,i,j} = \langle v_n,v_{i}v_{j} \rangle_{4\pi}, \quad \forall
(n,i,j) \in \mathbb{N}^3.
$$
According to Lemmas A.1 and A.2 of Appendix A of \cite{BSTU}, if
$W \in L^2_{\rm per}([-2\pi,2\pi])$ and $\| v_n \|^2_{L^2_{\rm per}([-2\pi,2\pi])}
= 1$ for any $n \in \mathbb{N}$, then there exists an $n$-independent constant
$C > 0$, such that
\begin{equation}
\label{bound-K} |K_{n,i,j}| \leq \frac{C}{(1 +
|n-i-j|)^p}, \qquad \forall (n,i,j) \in \mathbb{N}^3,
\end{equation}
for any $0 < p < 2$. By explicit computation, we obtain
\begin{eqnarray*}
\| \vec{\mbox{\boldmath $\phi$}} \star \vec{\mbox{\boldmath $\varphi$}}
\|_{l^1_s(\mathbb{N})} & \leq & \sum_{n \in \mathbb{N}} (1 + n)^s
\sum_{i \in \mathbb{N}} \sum_{j \in \mathbb{N}}
|K_{n,i,j}| |\phi_{i}| |\varphi_{j}| \\
& \leq & C \sum_{i \in \mathbb{N}} (1 + i)^s |\phi_{i}|
\sum_{j \in \mathbb{N}} (1 + j)^s |\varphi_{j}| \sum_{n \in
\mathbb{N}} \left( \frac{1 + n}{(1 + i)(1 +
j)} \right)^s \frac{1}{(1 + |n-i-j|)^p} \\
& \leq & \tilde{C} \sum_{i \in \mathbb{N}} (1 + i)^s
|\phi_{i}| \sum_{j \in \mathbb{N}} (1 + j)^s |\varphi_{j}|
\sum_{n \in \mathbb{N}} \left( 1 + \frac{n^s}{(1 + i)^s (1 +
j)^s} \right) \frac{1}{(1 + n)^p}
\end{eqnarray*}
for some $\tilde{C} > 0$. If $p > 1$ and $p - s > 1$, the bound is completed as
follows
$$
\| \vec{\mbox{\boldmath $\phi$}} \star \vec{\mbox{\boldmath $\varphi$}}
\|_{l^1_s(\mathbb{N})} \leq \tilde{C}_1 \| \vec{\mbox{\boldmath
$\phi$}}\|_{l^1_s(\mathbb{N}) } \|\vec{\mbox{\boldmath
$\varphi$}}\|_{l^1_s(\mathbb{N})} + \tilde{C}_2 \| \vec{\mbox{\boldmath
$\phi$}}\|_{l^1(\mathbb{N}) } \|\vec{\mbox{\boldmath
$\varphi$}}\|_{l^1(\mathbb{N})} \leq \left( \tilde{C}_1 +
\tilde{C}_2 \right) \| \vec{\mbox{\boldmath
$\phi$}}\|_{l^1_s(\mathbb{N}) } \|\vec{\mbox{\boldmath
$\varphi$}}\|_{l^1_s(\mathbb{N})}
$$
for some $\tilde{C}_1,\tilde{C}_2 > 0$. Since $1 < p < 2$, the
parameter $s$ must satisfy $0 < s < 1$.
\end{proof}

By using this construction, we prove the first main result of our
analysis.

\begin{theorem}
Let $W(x)$ be a bounded, even and $2\pi$-periodic function on $x
\in \mathbb{R}$. The nonlinear elliptic problem
(\ref{bifurcation}) has a continuous and $4 \pi$-periodic solution
$\Phi(x)$ for sufficiently small $\epsilon$ if there exists a
solution for $(A_1,A_2,A_3) \in \mathbb{C}^3$ of the algebraic
coupled-mode equations
\begin{equation}
\label{algebraic-system-cme}
\left\{ \begin{array}{rcl} (\Omega - \beta_1) A_1 & = & \sigma
\left[ \gamma_1 |A_1|^2 A_1 + \gamma_2 (2 |A_2|^2 A_1 + A_2^2
\bar{A}_1) + \gamma_3 (2
|A_3|^2 A_1 + A_3^2 \bar{A}_1) \right] + \epsilon R_1(A_1,A_2,A_3), \\
(\Omega - \beta_1) A_2 & = & \sigma \left[ \gamma_1 |A_2|^2 A_2 +
\gamma_2 (2 |A_1|^2 A_2 + A_1^2 \bar{A}_2) + \gamma_3 (2
|A_3|^2 A_2 + A_3^2 \bar{A}_2) \right]  + \epsilon R_2(A_1,A_2,A_3), \\
(\Omega - \beta_2) A_3 & = & \sigma \left[ \gamma_4 |A_3|^2 A_3 +
2 \gamma_3 (|A_1|^2 + |A_2|^2) A_3 + \gamma_3 (A_1^2 + A_2^2)
\bar{A}_3 \right] + \epsilon R_3(A_1,A_2,A_3),   \end{array}
\right.
\end{equation}
where parameters are given by
$$
\beta_1 = \langle\psi_1, W \psi_1\rangle_{4\pi} +
\langle\varphi_2, W \varphi_2\rangle_{4\pi}, \quad \beta_2 = 2
\langle\varphi_1, W \varphi_1\rangle_{4\pi}
$$
and
$$
\gamma_1 =
\langle\psi_1^2,\psi_1^2\rangle_{4\pi}\langle\varphi_2^2,\varphi_2^2\rangle_{4\pi},
\;\; \gamma_2 = |\langle\psi_1^2,\varphi_2^2\rangle_{4\pi}|^2,
\;\; \gamma_3 = \langle\psi_1^2,\varphi_1^2\rangle_{4\pi}
\langle\varphi_1^2,\varphi_2^2\rangle_{4\pi}, \;\; \gamma_4 =
|\langle\varphi_1^2,\varphi_1^2\rangle_{4\pi}|^2
$$
and the residual terms $R_{1,2,3}(A_1,A_2,A_3)$ are analytic
functions of $\epsilon$ near $\epsilon = 0$ satisfying the bounds
$$
\forall |\epsilon| < \epsilon_0 : \quad |R_{1,2,3}(A_1,A_2,A_3)|
\leq C_{1,2,3}( |A_1| + |A_2| + |A_3|),
$$
for some constants $C_{1,2,3} > 0$ and sufficiently small
$\epsilon_0$. Moreover, there exists an $\epsilon$-independent
constant $C > 0$ such that
$$
\| \Phi - A_1 \Phi_1 - A_2 \Phi_2 - A_3 \Phi_3
\|_{C^0_b(\mathbb{R}^2)} \leq C \epsilon,
$$
where $(\Phi_1,\Phi_2,\Phi_3)$ are the modes defined by Lemma
\ref{lemma-resonance}. \label{theorem-1}
\end{theorem}

\begin{proof}
$4\pi-$periodic solutions of the nonlinear elliptic problem (\ref{bifurcation})
are expanded in the form
\begin{equation}
\label{decomposition-2} \Phi(x) = \sum_{(n_1,n_2) \in
\mathbb{N}^2} \Phi_{n_1,n_2} v_{n_1}(x_1) v_{n_2}(x_2).
\end{equation}
Let the vector $\vec{\mbox{\boldmath $\Phi$}}$ with the elements
of the sequence $\{ \Phi_{n_1,n_2} \}_{(n_1,n_2) \in \mathbb{N}^2}$
belong to the vector space $l^1_s(\mathbb{N}^2)$
equipped with the norm
\begin{equation}
\label{vector-space-1-2d} \| \vec{\mbox{\boldmath
$\Phi$}}\|_{l^1_s(\mathbb{N}^2)} = \sum_{(n_1,n_2) \in
\mathbb{N}^2} (1 + n_1)^s (1+n_2)^s |\Phi_{n_1,n_2}| < \infty.
\end{equation}
Similarly to Lemma \ref{corollary-decomposition-1}, it follows
that the function $\Phi(x)$ is continuous if $\vec{\mbox{\boldmath
$\Phi$}} \in l^1_s(\mathbb{N}^2)$ with $s > \frac{1}{2}$. Indeed,
\begin{eqnarray*}
\| \Phi \|_{C^0_b(\mathbb{R}^2)} & \leq & \sum_{(n_1,n_2) \in
\mathbb{N}^2} |\Phi_{n_1,n_2}| \| v_{n_1} \|_{C^0_b([-2\pi,2\pi])}
\| v_{n_2} \|_{C^0_b([-2\pi,2\pi])} \\ & \leq & C \sum_{(n_1,n_2)
\in \mathbb{N}^2} |\Phi_{n_1,n_2}| \| v_{n_1}
\|_{H^s_{\rm per}([-2\pi,2\pi])} \| v_{n_2} \|_{H^s_{\rm per}([-2\pi,2\pi])}
\\ & \leq & \tilde{C} \sum_{(n_1,n_2) \in
\mathbb{N}^2} (1 + n_1)^s (1 + n_2)^s |\Phi_{n_1,n_2}|
\end{eqnarray*}
for some $C, \tilde{C} > 0$ and any $s > \frac{1}{2}$. By
construction, the solution $\Phi(x)$ is $4\pi$-periodic in both
coordinates of $x \in \mathbb{R}^2$. The partial differential
equation (\ref{bifurcation}) is rewritten in the lattice form,
which is diagonal with respect to the linear terms
\begin{equation}
\label{bifurcation-lattice} \left[ \nu_{n_1} + \nu_{n_2} -
\omega_0 - \epsilon \Omega \right] \Phi_n = - \epsilon \sigma
\sum_{(m,i,j) \in \mathbb{N}^6} M_{n,m,i,j} \Phi_m \bar{\Phi}_i
\Phi_j,  \quad \forall n = (n_1,n_2) \in \mathbb{N}^2,
\end{equation}
where $\nu_n$ depend on $\eta = \eta_0 + \epsilon$ and
$$
M_{n,m,i,j} = \langle v_{n_1} v_{i_1}, v_{m_1} v_{j_1}
\rangle_{4\pi} \; \langle v_{n_2} v_{i_2}, v_{m_2} v_{j_2}
\rangle_{4\pi}, \quad \forall (n_1,n_2,m_1,m_2,i_1,i_2,j_1,j_2)
\in \mathbb{N}^8.
$$
Since $M_{n,m,i,j}$ is a product of one-dimensional inner products
and the weights in the norm in $l^1_s(\mathbb{N}^2)$ are
separable, Lemma \ref{lemma-vector-field} applies and guarantees
that the nonlinear vector field of the lattice equations
(\ref{bifurcation-lattice}) is closed in $l^1_s(\mathbb{N}^2)$ for
$0 < s < 1$. Therefore, solutions of the lattice equations
(\ref{bifurcation-lattice}) can be considered in the space
$l^1_s(\mathbb{N}^2)$ for any $\frac{1}{2} < s < 1$.

By Lemma \ref{lemma-resonance}, the three resonant modes
are isolated from all other $4\pi$-periodic modes. By Lemma
\ref{lemma-perturbation-theory}, the eigenvalues $\nu_{n_1} +
\nu_{n_2}$ are then analytic in $\eta$ near $\eta = \eta_0$. The
three resonant modes correspond to the values of $n$ in the set
$\{ (1,3) ; (3,1) ; (2,2) \}$. Therefore, we decompose the
solution $\vec{\mbox{\boldmath $\Phi$}}$ of the lattice equations
(\ref{bifurcation-lattice}) into
\begin{equation}
\label{decomposition-3} \vec{\mbox{\boldmath $\Phi$}} = A_1 {\bf
e}_{1,3} + A_2 {\bf e}_{3,1} + A_3 {\bf e}_{2,2} +
\vec{\mbox{\boldmath $\Psi$}},
\end{equation}
where $\{ {\bf e}_{1,3},{\bf e}_{3,1},{\bf e}_{2,2} \}$ are unit
vectors on $\mathbb{N}^2$, ${\rm Span}({\bf e}_{1,3},{\bf
e}_{3,1},{\bf e}_{2,2})$ is the kernel of the linearized system at
the zero solution for $\epsilon = 0$, and $\vec{\mbox{\boldmath
$\Psi$}}$ lies in the orthogonal complement of the kernel such
that $\Psi_{1,3} = \Psi_{3,1} = \Psi_{2,2} = 0$. The linearized
operator projected onto the orthogonal complement of the kernel is
continuously invertible for sufficiently small $\epsilon$. By the
Implicit Function Theorem in the space $l^1_s(\mathbb{N}^2)$ for
any $\frac{1}{2} < s < 1$, there exists a unique map
$\vec{\mbox{\boldmath $\Psi$}}_{\epsilon}(A_1,A_2,A_3) :
\mathbb{C}^3 \mapsto l^1_s(\mathbb{N}^2)$ for sufficiently small
$\epsilon$. Moreover, the map is locally analytic in $\epsilon$
near $\epsilon = 0$, such that $\vec{\mbox{\boldmath
$\Psi$}}_0(A_1,A_2,A_3) = {\bf 0}$. In addition,
$\vec{\mbox{\boldmath $\Psi$}}_{\epsilon}(0,0,0) = {\bf 0}$ for
any $\epsilon \in \mathbb{R}$ and $(\partial_{\epsilon}
\vec{\mbox{\boldmath
$\Psi$}}_\epsilon(A_1,A_2,A_3))|_{\epsilon=0}$ is a homogeneous
cubic polynomial of $(A_1,A_2,A_3)$. Let $\delta \equiv |A_1| +
|A_2| + |A_3| < \delta_0$ for a fixed $\epsilon$-independent
$\delta_0 > 0$. Then the map $\vec{\mbox{\boldmath
$\Psi$}}_{\epsilon}(A_1,A_2,A_3)$ satisfies the bound
\begin{equation}
\label{bound-remainder} \forall |\epsilon| < \epsilon_0 : \quad \|
\vec{\mbox{\boldmath $\Psi$}}_{\epsilon}(A_1,A_2,A_3)
\|_{l^1_s(\mathbb{N}^2)} \leq \epsilon C \left( |A_1| + |A_2| +
|A_3| \right),
\end{equation}
where $\epsilon_0 > 0$ is sufficiently small, $\delta_0 > 0$ is
finite, and the constant $C > 0$ is independent of $\epsilon$ and
$\delta$.

We can now consider the three equations of the system
(\ref{bifurcation-lattice}) for $n = \{ (1,3) ; (3,1) ; (2,2) \}$.
After the Taylor expansion of (\ref{bifurcation-lattice}) in
$\epsilon$ the formula (\ref{leading-order}) yields the
coefficients $\beta_1,\beta_2$ since
$\partial_\eta=\partial_\epsilon$. Note that the $2\pi$-inner
product in (\ref{leading-order}) is replaced by the $4\pi$-inner
product. Although the Bloch modes $v_n(x)$ are now normalized over the $4\pi$-long
interval, the value of (\ref{leading-order}) is unchanged. Using
the bound (\ref{bound-remainder}) and the fact that the lattice
equations are closed in $l^1_s(\mathbb{N}^2)$, we derive the
algebraic coupled-mode equations in Theorem \ref{theorem-1}. The
coefficients of these equations can be simplified due to the fact
that $W(x)$ is even on $x \in \mathbb{R}$, such that the
eigenfunctions $\psi_1(x)$ and $\varphi_1(x)$ are even while
$\varphi_2(x)$ is odd on $x \in \mathbb{R}$. As a result, many
coefficients of the coupled-mode system are identically zero, for instance
$\langle \psi_1 \varphi_2,\varphi_1^2 \rangle_{4\pi} = 0$, $\langle
\varphi_1^2,\varphi_1 \varphi_2 \rangle_{4\pi} = 0$ and so on. The
residual terms $R_{1,2,3}(A_1,A_2,A_3)$ are estimated from the map
$\vec{\mbox{\boldmath $\Psi$}}_{\epsilon}(A_1,A_2,A_3)$ with the
bound (\ref{bound-remainder}). The last estimate in the statement of Theorem
\ref{theorem-1} is obtained from (\ref{bound-remainder}) and Lemma
\ref{corollary-decomposition-1}.
\end{proof}

We shall refer to the system (\ref{algebraic-system-cme}) without remainder terms
$\epsilon R_{1,2,3}(A_1,A_2,A_3)$ as to the truncated coupled-mode system.

\begin{corollary}
There exist five invariant reductions of the truncated coupled-mode system,
namely (i) $A_1 = A_2 = 0$, (ii) $A_1
= A_3 = 0$, (iii) $A_2 = A_3 = 0$, (iv) $A_1 = 0$, and (v) $A_2 =
0$, which persist in the full lattice equations
(\ref{bifurcation-lattice}).
\end{corollary}

\begin{proof}
Any of the five invariant reductions implies symmetry constraints
on the function $\Phi(x)$ at the leading order of the
decomposition (\ref{decomposition-3}). For instance, if $A_1 = A_2
= 0$, then $\Phi(x)$ is $2\pi$-antiperiodic with respect to both
$x_1$ and $x_2$; if $A_1 = 0$, then $\Phi(x)$ is
$2\pi$-antiperiodic in $x_1$ and $4 \pi$-periodic in $x_2$, and so
on. By the completeness results of Proposition
\ref{proposition-1}, all other terms of the decomposition
(\ref{decomposition-2}) which violate the symmetry constraints on
the solution $\Phi(x)$ can be set to be identically zero. By the
Implicit Function Theorem, the zero solution is unique near
$\epsilon = 0$. Therefore, the series (\ref{decomposition-2})
shrinks to fewer terms, the reduction persists for sufficiently
small $\epsilon$, and the proof of Theorem \ref{theorem-1}
applies.
\end{proof}

\begin{remark}
{\rm The invariant reduction $A_3 = 0$ of the truncated coupled-mode system
may not satisfy the full lattice equations
(\ref{bifurcation-lattice}) since the solution $\Phi(x)$ with $A_3
= 0$ is still a $4 \pi$-periodic function of both $x_1$ and $x_2$
and the series (\ref{decomposition-2}) can not be shrunk to fewer terms
for $A_3 = 0$.}
\end{remark}

\begin{remark}
\label{remark-l2} {\rm The Fourier--Bloch decomposition needed in
Theorem \ref{theorem-1} can be alternatively developed for
$\vec{\mbox{\boldmath $\phi$}}$ in the vector space $l^2_s(\mathbb{N})$
equipped with the squared norm
\begin{equation}
\label{vector-space-1a} \| \vec{\mbox{\boldmath $\phi$}}
\|^2_{l^2_s(\mathbb{N})} = \sum_{n \in \mathbb{N}} (1+n^2)^s
|\phi_n|^2 < \infty.
\end{equation}
Indeed, the squared norm $\| \phi \|^2_{H^s(\mathbb{P}_k)}$ is
equivalent to the integral
\begin{eqnarray*}
\int_{0}^{2\pi} \left| (c_{\eta} + L_{1D})^{s/2} \phi(x)
\right|^2 dx & = & \sum_{n_1 \in \mathbb{N}} \sum_{n_2 \in
\mathbb{N}} \phi_{n_1} \bar{\phi}_{n_2} \left( c_{\eta} +
\rho_n(k) \right)^{s} \langle
u_{n_1}(\cdot;k),u_{n_2}(\cdot;k)\rangle_{2\pi} \\ & = & \sum_{n
\in \mathbb{N}} \left( c_{\eta} + \rho_n(k) \right)^{s}
|\phi_n|^2,
\end{eqnarray*}
where we have used the orthogonality relation
(\ref{orthogonality-1}). By the asymptotic distribution
(\ref{asymptotic-distribution}), $\| \phi\|^2_{H^s(\mathbb{P}_k)}$
is thus equivalent to $\| \vec{\mbox{\boldmath
$\phi$}}\|^2_{l^2_s(\mathbb{N})}$. By Sobolev's Embedding Theorem,
$\| \phi \|_{C_b^0(\mathbb{P}_k)} \leq C \| \phi
\|_{H^s(\mathbb{P}_k)}$ for some $C > 0$ and any $s >
\frac{1}{2}$. Therefore, the decomposition (\ref{decomposition-1})
produces a continuous function $\phi(x)$ on $x \in \mathbb{R}$
if $\vec{\mbox{\boldmath $\phi$}} \in l^2_s(\mathbb{N})$ with $s >
\frac{1}{2}$. Furthermore, using Lemma 3.4 in \cite{BSTU}, one can
prove that the nonlinear term maps $l^2_s(\mathbb{N})$ to
$l^2_s(\mathbb{N})$ for any $s > \frac{1}{2}$, such that the
arguments of the Lyapunov--Schmidt reductions of Theorem
\ref{theorem-1} will work in the space $l^2_s(\mathbb{N})$ for any
$s > \frac{1}{2}$. Note that this approach would lift the upper
bound on the index $s$ in Theorem \ref{theorem-1}, where
$\frac{1}{2} < s < 1$. }
\end{remark}

\begin{example} \label{example-2}
{\rm For $W(x) = 1 - \cos x$ and $\eta = \eta_0 \approx 0.1745$,
the parameters of the algebraic coupled-mode equations (\ref{algebraic-system-cme}) are
approximated numerically as follows:
$$
\beta_1 \approx 2.2835, \quad \beta_2 \approx 0.9183
$$
and
$$
\gamma_1 \approx 9.4829 \times 10^{-3}, \quad \gamma_2 \approx
4.5196 \times 10^{-3}, \quad \gamma_3 \approx 3.7942 \times
10^{-3}, \quad \gamma_4 \approx 1.5981 \times 10^{-2}.
$$
}
\end{example}

The algebraic coupled-mode equations describe both linear and
nonlinear corrections to the eigenvalues
$\omega = \rho_{n_1}(k_1)+\rho_{n_2}(k_2)$ for fixed values of $(k_1,k_2)$
at the vertex points $X$, $X'$ and $M$. In particular, these
corrections show how the band edges of the three resonant bands
split for $\epsilon \neq 0$ and deform due to nonlinear
interactions between resonant Bloch modes. We skip further
analysis of the algebraic coupled-mode equations and proceed with
the decaying solutions of the nonlinear elliptic problem
(\ref{bifurcation}) described by the differential coupled-mode
equations.

\section{Differential coupled-mode equations}\label{S:diff_CME}

We consider here bounded and decaying solutions $\Phi(x)$ of
the bifurcation problem (\ref{bifurcation}) for small $\epsilon$.
The decomposition of the decaying solution $\Phi(x)$ depends now on
the completeness of the Bloch modes of the one-dimensional
Sturm--Liouville problem (\ref{E:Bloch_eval_problem}) in
$L^2(\mathbb{R})$, where all Bloch modes for all $k \in
\mathbb{T}$ must be incorporated in the Fourier--Bloch
decomposition.

\begin{proposition}
Let $W(x)$ be a bounded and $2\pi$-periodic function. Let $\{
\rho_n(k) \}_{n \in \mathbb{N}}$ and $\{ u_n(x;k) \}_{n \in
\mathbb{N}}$ be sets of eigenvalues and eigenfunctions of the
Sturm--Liouville problem (\ref{E:Bloch_eval_problem}) on $k \in
\mathbb{T}$, such that
\begin{equation}
\label{orthogonality-2} \langle u_n(\cdot;k), u_{n'}(\cdot;k')
\rangle_{\mathbb{R}} = \delta_{n,n'} \delta(k - k'),  \qquad
\forall n,n' \in \mathbb{N}, \;\; \forall k, k' \in \mathbb{T},
\end{equation}
where $\langle f,g \rangle_{\mathbb{R}}= \int_{\mathbb{R}}
\bar{f}(x) g(x) dx$ and $\delta(x)$ is the Dirac's delta function
in the distribution sense. Then, there exists a unitary
transformation ${\cal T} : L^2(\mathbb{R}) \mapsto
l^2(\mathbb{N},L^2(\mathbb{T}))$ given by
\begin{equation}
\label{direct-integral} \forall \phi \in L^2(\mathbb{R}) : \quad
\hat{\phi}(k) = {\cal T} \phi, \quad \hat{\phi}_n(k) =
\int_{\mathbb{R}} \bar{u}_n(y;k) \phi(y) dy, \quad \forall n \in
\mathbb{N}, \; \forall k \in \mathbb{T},
\end{equation}
where $l^2(\mathbb{N},L^2(\mathbb{T}))$ is equipped with the squared norm
$$
\| \hat{\phi} \|_{l^2(\mathbb{N},L^2(\mathbb{T}))}^2 = \sum_{n \in \mathbb{N}} \int_{\mathbb{T}}
|\hat{\phi}_n(k)|^2 dk.
$$
The inverse transformation is given by
\begin{equation}
\label{inverse-integral} \forall \hat{\phi} \in
l^2(\mathbb{N},L^2(\mathbb{T})) : \quad \phi(x) = {\cal T}^{-1}
\hat{\phi} = \sum_{n \in \mathbb{N}} \int_{\mathbb{T}}
\hat{\phi}_n(k) u_n(x;k) dk, \quad \forall x \in \mathbb{R}.
\end{equation}
\label{proposition-2}
\end{proposition}

\begin{proof}
The original proof of the proposition can be found in \cite{gelfand}.
Orthogonality and completeness of the Bloch wave
functions $\{ u_n(x;k) \}_{n \in \mathbb{N}}$ on $k \in
\mathbb{T}$ is summarized in Theorem XIII.98 on p.304 in \cite{RS}.
Existence of the unitary transformation (\ref{direct-integral})
with the inverse transformation (\ref{inverse-integral}) is summarized in
Theorems XIII.97 and XIII.98 on pp. 303--304 in \cite{RS}.
\end{proof}

\begin{lemma}
\label{corollary-decomposition-2} Let $\phi(x)$ be defined by the
decomposition (\ref{inverse-integral}), where $\hat{\phi}$ belongs
to the vector space $l^1_s(\mathbb{N},L^1(\mathbb{T}))$ with the norm
\begin{equation}
\label{vector-space-2} \| \hat{\phi}
\|_{l^1_s(\mathbb{N},L^1(\mathbb{T}))} = \sum_{n \in \mathbb{N}}
(1+n)^s \| \hat{\phi}_n \|_{L^1(\mathbb{T})} = \sum_{n \in
\mathbb{N}} (1+n)^s \int_{\mathbb{T}} |\hat{\phi}_n(k)| dk  <
\infty.
\end{equation}
If $s > \frac{1}{2}$, then $\phi(x)$ is a continuous function on
$x \in \mathbb{R}$ such that $\phi(x) \to 0$ as $|x| \to \infty$.
\end{lemma}

\begin{proof}
The proof is similar to that of Lemma
\ref{corollary-decomposition-1} since the asymptotic bound
(\ref{asymptotic-distribution}) is uniform on $k \in \mathbb{T}$.
Therefore,
\begin{eqnarray*}
\| \phi \|_{C^0_b(\mathbb{R})} \leq \sum_{n \in \mathbb{N}}
\int_{\mathbb{T}} |\hat{\phi}_n(k)| \| u_n(\cdot;k)
\|_{C^0_b(\mathbb{P}_k)} dk & \leq & C \sum_{n \in \mathbb{N}}
\int_{\mathbb{T}} |\hat{\phi}_n(k)| \| u_n(\cdot;k)
\|_{H^s(\mathbb{P}_k)} dk \leq \tilde{C} \| \hat{\phi}
\|_{l^1_s(\mathbb{N},L^1(\mathbb{T}))},
\end{eqnarray*}
for some $C, \tilde{C} > 0$ and any $s > \frac{1}{2}$. The decay
of $\phi(x)$ as $|x| \to \infty$ follows from the Riemann--Lebesgue
Lemma applied to the Fourier--Bloch transform. Indeed, the summation in $n \in \mathbb{N}$
of the integrals on $k \in \mathbb{T}$ can be written as an integral
on $k \in \mathbb{R}$ in the form
$$
\phi(x) = \int_{\mathbb{R}} \hat{\phi}(k) v(x;k) d k,
$$
where $v(x;k)$ for $k \in \mathbb{R}$ is related to the Bloch functions 
$u_n(x;k)$ for $k \in \mathbb{T}$ (see \cite{BSTU} for details).
Since $v(x;k) = e^{ikx} w(x;k)$, where $w(x;k)$ is periodic in $x$ with period $2\pi$
and is uniformly bounded in $C^0_b(\mathbb{R})$ with respect to both $x$ and $k$,
the Riemann--Lebesque Lemma applies
to the Fourier--Bloch transform in the same manner as it applies to the standard Fourier
transform if $\hat{\phi} \in L^1(\mathbb{R})$.
\end{proof}

\begin{lemma}
Let $W(x)$ be a bounded, piecewise-continuous and $2\pi$-periodic
function. Let $\rho_n(k_0)$ be the extremal value of $\rho_n(k)$
for either $k_0 = 0$ or $k_0 = \pm \frac{1}{2}$ and assume that
the adjacent bands are bounded away from the value $\rho_n(k_0)$.
Then, $\rho_n(k)$ is an analytic function at the point $k =
k_0$, such that
\begin{equation}
\rho'(k_0) = 0, \quad \rho_n''(k_0) \neq 0, \quad \left|
\rho_n(k) - \rho_n(k_0) - \frac{1}{2} \rho_n''(k_0)
(k-k_0)^2 \right| \leq C (k-k_0)^4
\end{equation}
for some $C > 0$ uniformly in $k \in \mathbb{T}$.
\label{lemma-band-curvature}
\end{lemma}

\begin{proof}
By Theorem XIII.95 in \cite{RS}, the function $\rho_n(k)$
for an isolated spectral band can be extended onto a smooth
Riemann surface which has no singularities other than square root
branch points at $k = 0$ and $k = \pm \frac{1}{2}$.
Therefore, the function $\rho_n(k)$ is expanded
in even powers of $(k-k_0)$ for $k_0 = 0$ or $k_0 = \pm
\frac{1}{2}$ and $\rho_n''(k_0) \neq 0$.
\end{proof}

\begin{lemma}
\label{lemma-vector-field-2} Let $W(x)$ be a bounded and
$2\pi$-periodic function. Let $\hat{\phi} = {\cal T} \phi$ and
$\hat{g} = {\cal T} (|\phi|^2 \phi)$. If $0 < s < 1$, then there
exists a $C > 0$ such that $\| \hat{g}
\|_{l^1_s(\mathbb{N},L^1(\mathbb{T}))} \leq C \| \hat{\phi}
\|^3_{l^1_s(\mathbb{N},L^1(\mathbb{T}))}$.
\end{lemma}

\begin{proof}
The proof is similar to that of Lemma \ref{lemma-vector-field}
since the asymptotic bound (\ref{bound-K}) is uniform on $k \in
\mathbb{T}$ \cite{BSTU} and the length of $\mathbb{T}$ is finite. The Banach
algebra property (\ref{Banach-algebra}) holds in space
$l^1_s(\mathbb{N},L^1(\mathbb{T}))$ for any $0 < s < 1$. Indeed, let
$$
\left( \hat{\phi} \star \hat{\varphi} \right)_n(k) = \sum_{j_1 \in
\mathbb{N}} \sum_{j_2 \in \mathbb{N}} \int_{\mathbb{T}}
\int_{\mathbb{T}} K_{n,j_1,j_2}(k,k_1,k_2) \hat{\phi}_{j_1}(k_1)
\hat{\varphi}_{j_2}(k_2) dk_1 dk_2, \qquad n \in \mathbb{N}, \; k
\in \mathbb{T},
$$
where $K_{n,j_1,j_2}(k,k_1,k_2) = \langle
u_n(\cdot;k),u_{j_1}(\cdot;k_1) u_{j_2}(\cdot;k_2)
\rangle_{\mathbb{R}^1}$ for all $(n,j_1,j_2) \in \mathbb{N}^3$ and
$(k,k_1,k_2) \in \mathbb{T}^3$. By an explicit computation, we
obtain
\begin{eqnarray*}
\| \hat{\phi} \star \hat{\varphi}
\|_{l^1_s(\mathbb{N},L^1(\mathbb{T}))} & \leq & \sum_{n \in
\mathbb{N}} (1 + n)^s \sum_{j_1 \in \mathbb{N}} \sum_{j_2 \in
\mathbb{N}} \int_{\mathbb{T}} \int_{\mathbb{T}} \int_{\mathbb{T}}
|K_{n,j_1,j_2}(k,k_1,k_2)| |\hat{\phi}_{j_1}(k_1)| |\hat{\varphi}_{j_2}(k_2)| dk dk_1 dk_2\\
& \leq & C \sum_{j_1 \in \mathbb{N}} \sum_{j_2 \in \mathbb{N}} (1
+ j_1)^s (1 + j_2)^s  \int_{\mathbb{T}} \int_{\mathbb{T}}
|\hat{\phi}_{j_1}(k_1)| |\hat{\varphi}_{j_2}(k_2)|  dk_1 dk_2  \\
& \phantom{t} & \phantom{text} \times  \sum_{n \in \mathbb{N}}
\left( \frac{1 + n}{(1 + j_1)(1 + j_2)} \right)^s \frac{1}{(1 +
|n-j_1-j_2|)^p}
\end{eqnarray*}
for some $C > 0$ and any $0 < p < 2$. The last sum is finite if $p
> 1$ and $p-s > 1$ (that is $0 < s < 1$), similarly to the proof of Lemma
\ref{lemma-vector-field}.
\end{proof}

We shall now define the coupled-mode system, which is a
central element of our paper. Let $A_{1,2,3}(y)$ be functions of
$y = \sqrt{\epsilon} x$ on $x \in \mathbb{R}^2$ which satisfy the
differential coupled-mode equations
\begin{equation}\label{E:CME_diff}
\left\{ \begin{array}{lcl} (\Omega - \beta_1) A_1 + \left(
\alpha_1
\partial_{y_1}^2 + \alpha_2 \partial_{y_2}^2 \right) A_1 & = &
\sigma \left[ \gamma_1 |A_1|^2 A_1 + \gamma_2 (2 |A_2|^2 A_1 +
A_2^2 \bar{A}_1)
+ \gamma_3 (2 |A_3|^2 A_1 + A_3^2 \bar{A}_1) \right] \\
(\Omega - \beta_1) A_2 + \left( \alpha_2 \partial_{y_1}^2 +
\alpha_1 \partial_{y_2}^2 \right) A_2 & = & \sigma \left[ \gamma_1
|A_2|^2 A_2 + \gamma_2 (2 |A_1|^2 A_2 + A_1^2 \bar{A}_2)
+ \gamma_3 (2 |A_3|^2 A_2 + A_3^2 \bar{A}_2) \right] \\
(\Omega - \beta_2) A_3 + \alpha_3 \left( \partial_{y_1}^2 +
\partial_{y_2}^2 \right) A_3 & = & \sigma \left[ \gamma_4 |A_3|^2
A_3 + 2 \gamma_3 (|A_1|^2 + |A_2|^2) A_3 + \gamma_3 (A_1^2 +
A_2^2) \bar{A}_3 \right],
\end{array} \right.
\end{equation}
where parameters $\beta_{1,2}$ and $\gamma_{1,2,3,4}$ are the same as
in Theorem \ref{theorem-1}, while
$$
\alpha_1 = \frac{1}{2} \rho''_1(0), \quad \alpha_2 = \frac{1}{2}
\rho_2''\left(\frac{1}{2}\right), \quad \alpha_3 = \frac{1}{2}
\rho_1''\left( \frac{1}{2}\right).
$$
Since $\lambda_1=\rho_1(0)$,
$\mu_1=\rho_1\left(\frac{1}{2}\right)$, and
$\mu_2=\rho_2\left(\frac{1}{2}\right)$ are non-degenerate
eigenvalues in the ordering (\ref{ordering}), the values of
$\alpha_{1,2,3}$ are non-zero according to Lemma
\ref{lemma-band-curvature}.

\begin{example} \label{example-3}
{\rm For $W(x) = 1 - \cos x$ and $\eta = \eta_0 = 0.1745$,
the parameters of the differential coupled-mode equations are
approximated numerically as follows:
$$
\alpha_1 \approx 0.9422, \quad \alpha_2 \approx 6.7813, \quad
\alpha_3 \approx -4.7890.
$$
}
\end{example}

Let $\hat{A}_{1,2,3}(p)$ on $p \in \mathbb{R}^2$ denote the
Fourier transforms of the functions $A_{1,2,3}(y)$ on $y \in
\mathbb{R}^2$ such that
\begin{equation}
\label{solution-integral} A_{1,2,3}(y) = \frac{1}{2\pi}
\int_{\mathbb{R}^2} \hat{A}_{1,2,3}(p) e^{i p \cdot y} dp, \qquad
\hat{A}_{1,2,3}(p) = \frac{1}{2\pi} \int_{\mathbb{R}^2}
A_{1,2,3}(y) e^{- i p \cdot y} dy.
\end{equation}
We shall use the standard norm in $L^1_q(\mathbb{R}^2)$ for the
Fourier transforms $\hat{A}_{1,2,3}(p)$ with any $q \geq 0$:
\begin{equation}
\| \hat{A} \|_{L^1_q(\mathbb{R}^2)} = \int_{\mathbb{R}^2} (1 +
|p|^2)^{q/2} |\hat{A}(p)| dp.
\end{equation}
The differential coupled-mode equations \eqref{E:CME_diff} are converted to the
equivalent integral form
\begin{eqnarray}
\label{cme-integral-form} \left\{ \begin{array}{ccc} \left(
\Omega - \beta_1 - \alpha_1 p_1^2 - \alpha_2 p_2^2 \right)
\hat{A}_1(p) & = & \sigma \hat{N}_1[\hat{A}_1,\hat{A}_2,\hat{A}_3](p),\\
\left( \Omega - \beta_1 - \alpha_2 p_1^2 - \alpha_1 p_2^2 \right)
\hat{A}_2(p) & = & \sigma \hat{N}_2[\hat{A}_1,\hat{A}_2,\hat{A}_3](p),  \\
\left( \Omega - \beta_2 - \alpha_3 p_1^2 - \alpha_3 p_2^2 \right)
\hat{A}_3(p) & = & \sigma
\hat{N}_3[\hat{A}_1,\hat{A}_2,\hat{A}_3](p), \end{array} \right.
\end{eqnarray}
where $\hat{N}_{1,2,3}[\hat{A}_1,\hat{A}_2,\hat{A}_3](p)$ denote
the Fourier transform of the cubic nonlinear terms of the
differential coupled-mode equations. We now prove the second main
result of our analysis.

\begin{theorem}
Let $W(x)$ be a bounded, piecewise-continuous, even and
$2\pi$-periodic function on $x \in \mathbb{R}$. Let $\frac{1}{4} <
r < \frac{1}{2}$. The nonlinear elliptic problem
(\ref{bifurcation}) has a continuous and decaying solution
$\Phi(x)$ for sufficiently small $\epsilon > 0$ if there exists a
solution $(\hat{B}_1,\hat{B}_2,\hat{B}_3) \in
L^1(D_{\epsilon},\mathbb{C}^3)$, compactly supported on the disk
$D_{\epsilon} = \{ p \in \mathbb{R}^2 : \; |p| <
\epsilon^{r-\frac{1}{2}}\} \subset \mathbb{R}^2$, which satisfies
the extended coupled-mode equations in the integral form
\begin{eqnarray}
\label{cme-integral-form-extended} \left\{ \begin{array}{ccc}
\left( \Omega - \beta_1 - \alpha_1 p_1^2 - \alpha_2 p_2^2 \right)
\hat{B}_1(p) - \sigma \hat{Q}_1(p) & = & \epsilon^{\tilde{r}}
\hat{R}_1(p), \\
\left( \Omega - \beta_1 - \alpha_2 p_1^2 - \alpha_1 p_2^2 \right)
\hat{B}_2(p) - \sigma \hat{Q}_2(p) & = & \epsilon^{\tilde{r}}
\hat{R}_2(p), \\
\left( \Omega - \beta_2 - \alpha_3 p_1^2 - \alpha_3 p_2^2 \right)
\hat{B}_3(p) - \sigma \hat{Q}_3(p) & = & \epsilon^{\tilde{r}}
\hat{R}_3(p), \end{array} \right.
\end{eqnarray}
where $\tilde{r} = \min(4r - 1,1 - 2r)$, $\hat{Q}_{1,2,3}(p)$
denote the cubic nonlinear terms
$\hat{N}_{1,2,3}[\hat{B}_1,\hat{B}_2,\hat{B}_3](p)$ truncated on
$p \in D_{\epsilon}$, and $\hat{R}_{1,2,3}(p)$ are bounded by
\begin{equation}
\label{bound-remainder-terms} \forall 0 < \epsilon < \epsilon_0 :
\quad \| \hat{R}_{1,2,3} \|_{L^1(D_{\epsilon})} \leq C_{1,2,3}
\left( \| \hat{B}_1 \|_{L^1(D_{\epsilon})} + \| \hat{B}_2
\|_{L^1(D_{\epsilon})} + \| \hat{B}_3 \|_{L^1(D_{\epsilon})}
\right),
\end{equation}
for some constants $C_{1,2,3} > 0$ and sufficiently small
$\epsilon_0$. Moreover, there exists an $\epsilon$-independent
constant $C
> 0$ such that
\begin{equation}
\label{bound-2} \| \Phi - B_1 \Phi_1 - B_2 \Phi_2 - B_3 \Phi_3
\|_{C^0_b(\mathbb{R}^2)} \leq C \epsilon^{1-2r},
\end{equation}
where $\Phi_{1,2,3}(x)$ are modes defined by Lemma
\ref{lemma-resonance} and $B_{1,2,3}(y)$ are defined by the
Fourier transform (\ref{solution-integral}). \label{theorem-2}
\end{theorem}

\begin{proof}
Solutions of the nonlinear elliptic problem (\ref{bifurcation})
are expanded in the form
\begin{equation}
\label{decomposition-3-2D} \Phi(x) = \sum_{(n_1,n_2) \in
\mathbb{N}^2} \int_{\mathbb{T}^2} \hat{\Phi}_{n_1,n_2}(k_1,k_2)
u_{n_1}(x_1;k_1) u_{n_2}(x_2;k_2) dk_1 dk_2.
\end{equation}
Let $\hat{\Phi}$ denote the union of functions
$\hat{\Phi}_{n_1,n_2}(k_1,k_2)$ for all $(n_1,n_2) \in \mathbb{N}^2$ and
$(k_1,k_2) \in \mathbb{T}^2$ and consider the vector space
$l^1_s(\mathbb{N}^2,L^1(\mathbb{T}^2))$ equipped with the norm
\begin{equation}
\label{vector-space-2-2d} \| \hat{\Phi}
\|_{l^1_s(\mathbb{N}^2,L^1(\mathbb{T}^2))} = \sum_{(n_1,n_2) \in
\mathbb{N}^2} (1 + n_1)^s (1+n_2)^s \int_{\mathbb{T}}
\int_{\mathbb{T}} |\hat{\Phi}_{n_1,n_2}(k_1,k_2)| dk_1 dk_2  <
\infty.
\end{equation}
One can show similarly to the proof of Theorem \ref{theorem-1}
that the statements of Lemmas \ref{corollary-decomposition-2} and
\ref{lemma-vector-field-2} extend in two dimensions to the vector space with the norm
(\ref{vector-space-2-2d}). Therefore, we convert the partial
differential equation (\ref{bifurcation}) to the integral form,
which is diagonal with respect to the linear terms
\begin{eqnarray}
 \label{bifurcation-integral} & \phantom{t} & \left[ \rho_{n_1}(k_1) +
\rho_{n_2}(k_2) - \omega_0 - \epsilon \Omega \right]
\hat{\Phi}_n(k) \\
 \nonumber& = & - \epsilon \sigma
\sum_{(m,i,j) \in \mathbb{N}^6} \int_{\mathbb{T}^6}
M_{n,m,i,j}(k,l,\kappa,\lambda) \hat{\Phi}_m(l)
\bar{\hat{\Phi}}_i(\kappa) \hat{\Phi}_j(\lambda) dl d \kappa d
\lambda
\end{eqnarray}
for all $n = (n_1,n_2) \in \mathbb{N}^2$ and $k = (k_1,k_2) \in
\mathbb{T}^2$, where
$$
M_{n,m,i,j}(k,l,\kappa,\lambda) = \langle u_{n_1}(\cdot;k_1)
u_{i_1}(\cdot;\kappa_1), u_{m_1}(\cdot;l_1)
u_{j_1}(\cdot;\lambda_1) \rangle_{\mathbb{R}} \; \langle
u_{n_2}(\cdot;k_2) u_{i_2}(\cdot;\kappa_2), u_{m_2}(\cdot;l_2)
u_{j_2}(\cdot;\lambda_2) \rangle_{\mathbb{R}},
$$
for all $(n_1,n_2,m_1,m_2,i_1,i_2,j_1,j_2) \in \mathbb{N}^8$ and
$(k_1,k_2,l_1,l_2,\kappa_1,\kappa_2,\lambda_1,\lambda_2) \in
\mathbb{T}^8$. In view of Lemmas \ref{corollary-decomposition-2} and
\ref{lemma-vector-field-2} we consider solutions of the integral equations
(\ref{bifurcation-integral}) in the space
$l^1_s(\mathbb{N}^2,L^1(\mathbb{T}^2))$ for any $\frac{1}{2} < s <
1$.

When $\eta = \eta_0$, the multiplication operator
$\rho_{n_1}(k_1) + \rho_{n_2}(k_2) - \omega_0$ vanishes at the
points $X'$, $X$, and $M$, where the values of $k$ are given by
\begin{equation}
\label{resonant-points} \left\{
\left( 0,\frac{1}{2} \right); \left( \frac{1}{2},0 \right); \left( \frac{1}{2},\frac{1}{2}
\right) \right\} \subset \mathbb{T}^2,
\end{equation}
with $n \in \{ (1,3) ; (3,1) ; (2,2) \} \subset \mathbb{N}^2$ resp.
Therefore, we apply the method of Lyapunov--Schmidt reductions \cite{GS84} and
decompose the solution $\hat{\Phi}$ of the integral equations (\ref{bifurcation-integral}) into
\begin{equation}
\label{decomposition-4} \hat{\Phi}(k) = \hat{U}_1(k) \chi_{D_1}(k)
{\bf e}_{1,3} + \hat{U}_2(k) \chi_{D_2}(k)  {\bf e}_{3,1} + \hat{U}_3(k)
\chi_{D_3}(k) {\bf e}_{2,2} + \hat{\Psi}(k),
\end{equation}
where $\{ {\bf e}_{1,3},{\bf e}_{3,1},{\bf e}_{2,2} \}$ are unit vectors on
$\mathbb{N}^2$, $D_{1,2,3}$ are disks of radius $\epsilon^r$
centered at the points $k$ of the set (\ref{resonant-points}), see Figure \ref{F:k_sp_decomp}.
Since $k \in \mathbb{T}^2$, where $\mathbb{T}^2$ is the first Brillouin zone,
we map $D$ on $\mathbb{T}^2$ by periodic continuation and denote the resulting object by
$D \cap \mathbb{T}^2$. In the representation (\ref{decomposition-4}),
$\chi_D(k)$ is the characteristic function on $k \in \mathbb{T}^2$
($\chi_D(k) = 1$ if $k \in D \cap \mathbb{T}^2$ and $\chi_D(k)
= 0$ if $k \notin D \cap \mathbb{T}^2$), and $\hat{\Psi}(k)$ is zero identically on
$k \in D_{1,2,3}$ for the corresponding values of $n$:
$$
\hat{\Psi}_{(1,3)}=0 \;\; \mbox{on} \;\; D_1, \quad \hat{\Psi}_{(3,1)}=0 \;\; \mbox{on} \;\;
D_2, \quad \hat{\Psi}_{(2,2)}=0 \;\; \mbox{on} \;\; D_3.
$$

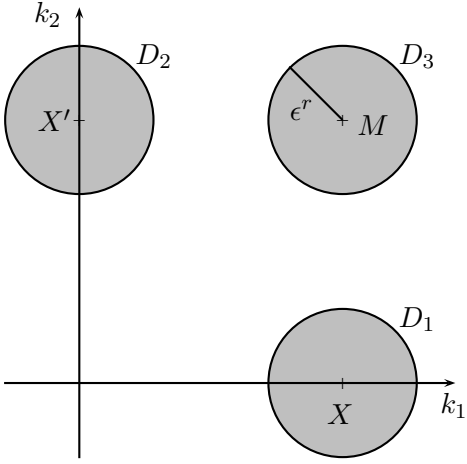
\begin{figure}[htpb]
\centering
\begin{minipage}[c]{.5\textwidth}
\centering
\setlength{\unitlength}{1cm}
\begin{picture}(6,6)(0,0)
  \pscircle[fillstyle=solid,fillcolor=lightgray](1,4.5){1}
  \pscircle[fillstyle=solid,fillcolor=lightgray](4.5,4.5){1}
  \pscircle[fillstyle=solid,fillcolor=lightgray](4.5,1){1}
  \psline{->}(0,1)(6,1)
  \psline{->}(1,0)(1,6)
  \psline(4.5,4.5)(3.79,5.21)
  \psdots*[dotstyle=+](1,4.5)(4.5,1)(4.5,4.5)
  \put(0.45,4.35){$X'$}
  \put(1.75,5.25){$D_2$}
  \put(4.3,0.45){$X$}
  \put(5.25,1.75){$D_1$}
  \put(4.7,4.3){$M$}
  \put(5.25,5.25){$D_3$}
  \put(3.8,4.5){$\epsilon^r$}
  \put(5.8,0.6){$k_1$}
  \put(0.4,5.8){$k_2$}
\end{picture}
 \end{minipage}
\begin{minipage}[c]{.45\textwidth}
\centering \caption{Decomposition of the $k$-space}
\label{F:k_sp_decomp}
\end{minipage}
\end{figure}

For any of the three resonant points of $k$, it follows by Lemma
\ref{lemma-band-curvature} that $\rho'_{n_j}(k_j) = 0$ and
$\rho''_{n_j}(k_j) \neq 0$ for $j = 1,2$ with the $n$ corrresponding to the resonance point $k$.
As a result, there exists a constant $C > 0$ such that
\begin{equation}
\label{lower-bound} \min_{(k_1,k_2) \in {\rm supp}(\hat{\Psi})
\atop (n_1,n_2) \in \mathbb{N}^2} |\rho_{n_1}(k_1) |_{\eta =
\eta_0} + \rho_{n_2}(k_2) |_{\eta = \eta_0} - \omega_0| \geq C
\epsilon^{2r}.
\end{equation}
Note that the minimal values of the multiplication operator occur
for the three resonant bands with the indices $n$ in the set $\{
(1,3) ; (3,1) ; (2,2) \}$. For all other non-resonant bands, the
multiplication operator is continuously invertible uniformly in
$\epsilon$ for sufficiently small $\epsilon$. If $2r < 1$, the
lower bound (\ref{lower-bound}) is still larger than the
perturbation terms of order $\epsilon$. By the Implicit Function
Theorem in the space $l^1_s(\mathbb{N}^2,L^1(\mathbb{T}^2))$ for
any $\frac{1}{2} < s < 1$, there exists a unique map
$\hat{\Psi}_{\epsilon}(\hat{U}_1,\hat{U}_2,\hat{U}_3) : L^1(D_1)
\times L^1(D_2) \times L^1(D_3) \mapsto
l^1_s(\mathbb{N}^2,L^1(\mathbb{T}^2))$ for sufficiently small
$\epsilon > 0$. The right-hand side of the integral equation
(\ref{bifurcation-integral}) is a homogeneous cubic polynomial of
$\hat{\Psi}(k)$ and $\hat{U}_1,\hat{U}_2,\hat{U}_3$ multiplied by
$\epsilon$. Let $\delta \equiv \|\hat{U}_1\|_{L^1(D_1)} +
\|\hat{U}_2 \|_{L^1(D_2)} + \|\hat{U}_3\|_{L^1(D_3)} < \delta_0$
for a fixed $\epsilon$-independent $\delta_0 > 0$. The map
$\hat{\Psi}_{\epsilon}(\hat{U}_1,\hat{U}_2,\hat{U}_3)$ satisfies
the bound
\begin{equation}
\label{bound-remainder-2d} \|
\hat{\Psi}_{\epsilon}(\hat{U}_1,\hat{U}_2,\hat{U}_3)
\|_{l^1_s(\mathbb{N}^2,L^1(\mathbb{T}^2))} \leq \epsilon^{1-2r} C
\left(\|\hat{U}_1\|_{L^1(D_1)} + \|\hat{U}_2 \|_{L^1(D_2)} +
\|\hat{U}_3\|_{L^1(D_3)} \right),
\end{equation}
for sufficiently small $0 < \epsilon < \epsilon_0$, finite $0 <
\delta < \delta_0$, and the $(\epsilon,\delta)$-independent
constant $C > 0$.

We shall now map the disks $D_{1,2,3}$ to the same disk
$D_{\epsilon} \subset \mathbb{R}^2$ in the stretched variable $p$
centered at the origin. To do so, we apply the scaling
transformation
\begin{equation}
\label{scaling-transformation} \hat{B}_j(p) = \epsilon
\hat{U}_j(k), \ p = \frac{k - k_0}{\epsilon^{1/2}}, \quad
\forall k \in D_j \cap \mathbb{T}^2, \quad j = 1,2,3,
\end{equation}
where $k_0$ is the corresponding point of the set
(\ref{resonant-points}). The new disk $D_{\epsilon}$ is now
$D_{\epsilon} = \{ p \in \mathbb{R}^2 : \; |p| <
\epsilon^{r-\frac{1}{2}}\}$ and it covers the entire plane $p \in
\mathbb{R}^2$ in the limit $\epsilon \to 0$ if $2r < 1$. Note that
the $L^1$-norm of $\hat{U}_j(k)$ on $k \in D_j \cap \mathbb{T}^2$
is invariant with
respect to the transformation (\ref{scaling-transformation}) in
the new variable $\hat{B}_j(p)$ on $p \in D_{\epsilon} \subset \mathbb{R}^2$,
such that $\|\hat{U}_j\|_{L^1(D_j)} = \| \hat{B}_j \|_{L^1(D_{\epsilon})}$
for any $j = 1,2,3$.

When the integral equations (\ref{bifurcation-integral})
are considered on the compact domains $D_1$, $D_2$, $D_3$ for the
three resonant modes and the scaling transformation
(\ref{scaling-transformation}) is applied to map these domains
into the domain $D_{\epsilon}$, we obtain the extended
coupled-mode system (\ref{cme-integral-form-extended}), where
remainder terms occur due to three different sources. The first
source comes from the component $\hat{\Psi} =
\hat{\Psi}_{\epsilon}(\hat{U}_1,\hat{U}_2,\hat{U}_3)$ and it is
estimated with the bound (\ref{bound-remainder-2d}). The remainder
terms have the order of $\epsilon^{1-2r}$ and are small if $2r <
1$. The second source comes from the expansion of $\rho_{n}(k)$
and $M_{n,m,i,j}(k,l,\kappa,\lambda)$ in powers of $\epsilon$
and it is estimated with Lemma \ref{lemma-perturbation-theory}. The
remainder terms have the order of $\epsilon^1$. The third source
comes from the expansion of $\rho_{n}(k)$ and
$M_{n,m,i,j}(k,l,\kappa,\lambda)$ in powers of $k - k_0$ and it is
estimated with Lemma \ref{lemma-band-curvature}. The remainder
terms have the order of $\epsilon^{4r-1}$ because of the following
estimate:
$$
\epsilon \| |p|^4 \hat{B}_j \|_{L^1(D_{\epsilon})} = \epsilon
\int_{D_{\epsilon}} |p|^4 \left| \hat{B}_j(p) \right| dp \leq
\epsilon^{4r - 1} \| \hat{B}_j \|_{L^1(D_{\epsilon})}, \quad j =
1,2,3.
$$
The remainder terms are small if $4r > 1$. Thus, the statement is
proved if $\frac{1}{4} < r < \frac{1}{2}$ and the largest
remainder terms have the order $\epsilon^{\tilde{r}}$ with
$\tilde{r} = \min(4 r - 1, 1 - 2 r)$.
\end{proof}

\begin{remark}
{\rm If $r = \frac{1}{3}$, then $\tilde{r} = r = \frac{1}{3}$ and
both remainder terms of the extended coupled-mode equations are of
the same order in $\epsilon$. The fastest convergence rate
of the remainder terms is $\epsilon^{1/3}$ because
$\max_{r\in (1/4,1/2)}(\tilde{r}) = 1/3$.}
\end{remark}

The differential coupled-mode system on $A_{1,2,3}(y)$ has
invariant reductions when one or two components of $A_{1,2,3}$ are
identically zero. In particular, the reduction $A_3 = 0$ recovers
coupled-mode equations near the band edge $C$ derived in
\cite{SY}, while the reduction $A_1 = A_2 = 0$ recovers
coupled-mode equations near the band edge $B$ in \cite{SY}.
However, these reductions do not generally persist in the extended
coupled-mode equations near $\eta = \eta_0$ since the reductions
do not imply any symmetry constraints on the function $\Phi(x)$.

The coupled-mode equations \eqref{cme-integral-form} linearized at the zero
solution describe the expansions of the function $\omega =
\rho_{n_1}(k_1) + \rho_{n_2}(k_2)$ near the vertex points
$X'$, $X$, and $M$ in terms of the perturbation parameter $\epsilon =
\eta - \eta_0$ and the small deviation of $k = (k_1,k_2)$ from the
value $k_0 = (k_{1,0},k_{2,0})$ in the set (\ref{resonant-points})
\begin{eqnarray}
\label{omega-expansion}
\omega = \omega_0 + \epsilon \partial_{\eta} \omega |_{\eta =
\eta_0,k=k_0} + \frac{1}{2} \rho_{n_1}''(k_{1,0}) (k_1 -
k_{1,0})^2 + \frac{1}{2} \rho_{n_2}''(k_{2,0}) (k_2 - k_{2,0})^2 +
\mathcal{O}(|k - k_0|^4).
\end{eqnarray}
We note that the values of $\partial_{\eta} \omega |_{\eta = \eta_0}$ are
equal for the vertex points $X$ and $X'$, while the values of
$\rho_{n_1}''(k_{1,0})$ and $\rho''_{n_2}(k_{2,0})$ have the same
sign for all points $X$, $X'$ and $M$. Moreover, $\text{sign}(\rho''_{n_j}(X_j))=\text{sign}(\rho''_{n_j}(X_j')), \ j=1,2$. A simple analysis of (\ref{omega-expansion}) shows that
the band gap between the three resonant Bloch band surfaces exists if
$$
\rho_{n_1}''(k_1) |_X > 0, \;\; \rho_{n_1}''(k_1) |_{M} < 0, \;\;\mbox{for}
\;\; \epsilon \partial_{\eta} \omega |_X > \epsilon
\partial_{\eta} \omega |_M
$$
or
$$
\rho_{n_1}''(k_1) |_X< 0, \;\; \rho_{n_1}''(k_1) |_M > 0, \;\;\mbox{for}
\;\; \epsilon \partial_{\eta} \omega |_X < \epsilon
\partial_{\eta} \omega |_M.
$$
Examples \ref{example-2} and \ref{example-3} show that the first
case occurs for the particular potential $W(x) = 1 - \cos x$ if
$\eta > \eta_0$, such that the band gap opens for $\epsilon > 0$ in the
interval $\beta_2 < \Omega < \beta_1$, in a correspondence with a
narrow band gap for $\eta > \eta_0$ in Figure \ref{F:bands}(b).

\section{Persistence of localized reversible solutions}\label{S:reversible_GS}

The coupled-mode system in the integral form
(\ref{cme-integral-form}) is different from the extended
coupled-mode system (\ref{cme-integral-form-extended}) in two
aspects. First, the convolution integrals are truncated in \eqref{cme-integral-form-extended} on the
domain $D_{\epsilon}$ for $|p| < \epsilon^{r - \frac{1}{2}}$,
where $\frac{1}{4} < r < \frac{1}{2}$. Second, the remainder terms
of order $\epsilon^{\tilde{r}}$ with $\tilde{r} = \min(4r - 1,1 -
2r)$ are present in \eqref{cme-integral-form-extended}. The first source is small if solutions of the
extended coupled-mode system are considered on $p \in
\mathbb{R}^2$ in the space $L^1_q(\mathbb{R}^2)$ for some $q \geq
0$. The second source is handled with the Implicit Function
Theorem. We shall consider here a special class of decaying solutions of the
coupled-mode system (\ref{cme-integral-form}) on $p \in
\mathbb{R}^2$ which lead to the corresponding solutions of the extended
coupled-mode system (\ref{cme-integral-form-extended}) on $p \in
D_{\epsilon}$.

\begin{definition}
\label{definition-reversible} A solution $(A_1,A_2,A_3)$ of the
differential coupled-mode system \eqref{E:CME_diff} on $y \in
\mathbb{R}^2$ is called reversible if it satisfies one
of the two constraints (\ref{reversibility}) and
(\ref{reversibility-2}) for each function $A_1(y)$, $A_2(y)$ and
$A_3(y)$.
\end{definition}

For notational simplicity, we write the differential coupled-mode
equations \eqref{E:CME_diff} in the form ${\bf F}({\bf A}) = {\bf 0}$, where ${\bf F}$
is a nonlinear operator on $(A_1,A_2,A_3) \in
C^2(\mathbb{R}^2,\mathbb{C}^3)$ with the range in
$C^0(\mathbb{R}^2,\mathbb{C}^3)$. The Jacobian of the operator
${\bf F}({\bf A})$ denoted by $D_{\bf A} {\bf F}({\bf A})$ is a
matrix differential operator, which is diagonal with respect to
the unbounded differential part and full with respect to the local
potential part. In many cases, the Jacobian operator $D_{\bf A}
{\bf F}({\bf A})$ can be block-diagonalized and simplified but we can
work with a general operator to prove the third main
result of our analysis.

\begin{theorem}
\label{theorem-3} Let $W(x)$ satisfy the same assumptions as in
Theorem \ref{theorem-2}. Let $\Omega$
belong to the interior of the band gap of the coupled-mode system
${\bf F}({\bf A}) = {\bf 0}$. Let $(A_1,A_2,A_3)$ on $y \in
\mathbb{R}^2$ be a reversible solution of the differential
coupled-mode system ${\bf F}({\bf A}) = {\bf 0}$ in the sense of
Definition \ref{definition-reversible} such that their Fourier
transforms satisfy $\hat{\bf A} \in
L^1_q(\mathbb{R}^2,\mathbb{C}^3)$ for all $q \geq 0$. Assume that the Jacobian operator
$D_{\bf A} {\bf F}({\bf A})$ has a three-dimensional kernel with
the eigenvectors $\{
\partial_{y_1} {\bf A}, \partial_{y_2} {\bf A}, i {\bf A} \}$.
Then, there exists a solution $(\hat{B}_1,\hat{B}_2,\hat{B}_3) \in
L^1(D_{\epsilon},\mathbb{C}^3)$ of the extended coupled-mode
system (\ref{cme-integral-form-extended}) such that
\begin{equation}
\label{bound-final} \forall 0 < \epsilon < \epsilon_0 : \quad \|
\hat{B}_j - \hat{A}_j \|_{L^1(D_{\epsilon})} \leq C_j
\epsilon^{\tilde{r}}, \quad \forall j = 1,2,3,
\end{equation}
for some $\epsilon$-independent constant $C_j > 0$ and
sufficiently small $\epsilon_0$.
\end{theorem}

\begin{proof}
First, we consider system (\ref{cme-integral-form-extended}) on
the entire plane $p \in \mathbb{R}^2$ in the space
$L^1_q(\mathbb{R}^2)$ for a fixed $q \geq 0$ by extending
the right-hand-side functions $\hat{R}_{1,2,3}(p)$ on $\mathbb{R}^2$
with a compact support on $D_{\epsilon}$. Thanks to the
assumption on the existence of a solution $\hat{\bf A}$ of
$\hat{\bf F}(\hat{\bf A}) = {\bf 0}$ and its fast decay in
$L^1_q(\mathbb{R}^2)$ for all $q \geq 0$, we decompose the
solution of the extended system $\hat{\bf F}(\hat{\bf B}) =
\epsilon^{\tilde{r}} \hat{\bf R}(\hat{\bf B})$ on $p \in
\mathbb{R}^2$ into $\hat{\bf B} = \hat{\bf A} + \hat{\bf b}$,
where $\hat{\bf b}$ solves the nonlinear problem in the form
\begin{equation}
\label{lin-abstract-form} \hat{J} \hat{\bf b} = \hat{\bf
N}(\hat{\bf b}), \quad \hat{J} = D_{\hat{\bf A}} \hat{\bf
F}(\hat{\bf A}), \;\; \hat{\bf N}(\hat{\bf b}) =
\epsilon^{\tilde{r}} \hat{\bf R}(\hat{\bf A} + \hat{\bf b}) -
\left[ \hat{\bf F}(\hat{\bf A} + \hat{\bf b}) - \hat{J} \hat{\bf
b} \right],
\end{equation}
where $\hat{J}$ is a linearized operator, $\hat{\bf F}(\hat{\bf
A} + \hat{\bf b}) - \hat{J} \hat{\bf b}$ is quadratic in $\hat{\bf
b}$, and $\hat{\bf R}(\hat{\bf A}+\hat{\bf b})$ maps an element $L_q^1(\mathbb{R}^2)$ to
itself. The kernel of the Jacobian operator $J = D_{\bf A} {\bf
F}({\bf A})$ is exactly three-dimensional, by the assumption, and
it is generated by the two-dimensional group of translations in $y
\in \mathbb{R}^2$ and a one-dimensional gauge invariance in
$\arg({\bf A})$. If $\Omega$ belongs to the interior of the band
gap of the coupled-mode system, the continuous spectrum of $D_{\bf
A} {\bf F}({\bf A})$ is bounded away from zero and the triple zero
eigenvalue is isolated from other eigenvalues of the discrete
spectrum of $D_{\bf A} {\bf F}({\bf A})$. The nonlinear elliptic
problem (\ref{stationary}) with a real-valued symmetric potential
$V(x)$ admits the gauge invariance $\phi(x) \to e^{i \alpha}
\phi(x)$ for any $\alpha \in \mathbb{R}$ and the reversibility
symmetries (\ref{reversibility}) and (\ref{reversibility-2}). The extended coupled-mode system
(\ref{cme-integral-form-extended}) obtained from the integral system
(\ref{bifurcation-integral}) after the
Lyapunov--Schmidt decomposition (\ref{decomposition-4})
and the scaling transformation (\ref{scaling-transformation}) inherits
all the symmetries, such that after restriction of ${\bf B}$ to functions satisfying (\ref{reversibility}) or
(\ref{reversibility-2}) and after setting uniquely the phase, e.g.
by ${\bf B}(0) \in \mathbb{R}$, system (\ref{lin-abstract-form}) is
formulated in the orthogonal complement of the kernel of $\hat{J}$. Therefore, the linearized
operator is continuously invertible for sufficiently small
$\epsilon > 0$ under this restriction. By the Implicit Function Theorem, there
exists a unique reversible solution in the form $\hat{\bf b} =
\hat{J}^{-1} \hat{\bf N}(\hat{\bf b})$ such that $\| \hat{\bf b}
\|_{L^1_q(\mathbb{R}^2,\mathbb{C}^3)} \leq C \epsilon^{\tilde{r}}$
for some $C > 0$ and any $q \geq 0$. This result implies the
desired bound (\ref{bound-final}) for $\hat{B}_j$ as a solution of
(\ref{cme-integral-form-extended}) on $p \in \mathbb{R}^2$.

Next, we consider the error between system
(\ref{cme-integral-form-extended}) on the disk $p \in
D_{\epsilon}$  and the same system on the plane $p \in
\mathbb{R}^2$. The error comes from the terms $\| \hat{\bf b}
\|_{L^1_{q+2}(D_{\epsilon}^{\perp},\mathbb{C}^3)}$ and $\|
\hat{\bf b} \|_{L^1_q(D_{\epsilon}^{\perp},\mathbb{C}^3)}$, where
$D_{\epsilon}^\perp = \mathbb{R}^2 \backslash D_{\epsilon}$, since
$\hat{\bf A} \in L^1_q(\mathbb{R}^2,\mathbb{C}^2)$ for all $q \geq
0$ and the cubic nonlinear terms $\hat{\bf N}_j(\hat{\bf b})$ of the
coupled-mode system (\ref{cme-integral-form}) map elements of
$L^1_q(\mathbb{R}^2)$ to elements of $L^1_q(\mathbb{R}^2)$.
Since $\hat{\bf b} = \hat{J}^{-1} \hat{\bf
N}(\hat{\bf b})$ and $\hat{J}^{-1}$ is a map from
$L^1_q(\mathbb{R}^2,\mathbb{C}^3)$ to
$L^1_{q+2}(\mathbb{R}^2,\mathbb{C}^3)$ for any $q \geq 0$, we
obtain that
$$
\forall \hat{\bf b} = \hat{J}^{-1} \hat{\bf N}(\hat{\bf b}) \in
L^1_q(\mathbb{R}^2,\mathbb{C}^3) : \quad \| \hat{\bf b}
\|_{L^1_{q+2}(D_{\epsilon}^{\perp},\mathbb{C}^3)} \leq \| \hat{\bf
b} \|_{L^1_{q+2}(\mathbb{R}^2,\mathbb{C}^3)}  \leq C \| \hat{\bf
N}(\hat{\bf b}) \|_{L^1_q(\mathbb{R}^2,\mathbb{C}^3)} \leq
\tilde{C} \epsilon^{\tilde{r}}
$$
for some $C,\tilde{C} > 0$. Since $\| \hat{\bf
b} \|_{L^1_{q}(D_{\epsilon}^{\perp},\mathbb{C}^3)}  \leq \| \hat{\bf
b} \|_{L^1_{q}(\mathbb{R}^2,\mathbb{C}^3)}  \leq C \epsilon^{\tilde{r}}$,
the error between system
(\ref{cme-integral-form-extended}) on the disk $p \in
D_{\epsilon}$  and the same system on the plane $p \in
\mathbb{R}^2$ is of the same order as the right-hand-side terms of
system (\ref{cme-integral-form-extended}), such that the desired
bound (\ref{bound-final}) holds on $p \in D_{\epsilon}$.
\end{proof}

\begin{corollary}\label{C:reversibility}
The solution $\Phi(x)$ constructed in Theorems \ref{theorem-2} and
\ref{theorem-3} is a reversible localized solution of the
bifurcation problem (\ref{bifurcation}) satisfying one of the two
constraints (\ref{reversibility}) and (\ref{reversibility-2}) on
$x \in \mathbb{R}^2$.
\end{corollary}

\begin{proof}
Since the modes $\Phi_{1,2,3}(x)$ of Lemma \ref{lemma-resonance}
and the solution $B_{1,2,3}(y)$ of Theorem \ref{theorem-3} satisfy
one of the two reversibility constraints (\ref{reversibility}) and
(\ref{reversibility-2}), the leading-order part of the
representation (\ref{decomposition-4}) produces a solution
$\Phi(x)$ which satisfies the same constraint. Since the symmetry
is also preserved in the integral equation
(\ref{bifurcation-integral}), the map
$\hat{\Psi}_{\epsilon}(\hat{U}_1,\hat{U}_2,\hat{U}_3)$ constructed
in Theorem \ref{theorem-2} inherits the symmetry and produces the
remainder term in the decomposition (\ref{decomposition-4})
satisfying the same reversibility constraint.
\end{proof}

There are several classes of localized reversible solutions of the
differential coupled-mode equations which satisfy the assumptions
of Theorem \ref{theorem-3}. These solutions are numerically
approximated in Section \ref{S:numerical}. According to Theorem
\ref{theorem-3}, all these solutions persist as localized
reversible solutions of the nonlinear elliptic problem
(\ref{stationary}).

\begin{corollary}\label{C:conv_rate_final}
The solution $\Phi(x)$ constructed in Theorems \ref{theorem-2} and
\ref{theorem-3} satisfy the bound:
\begin{equation}
\label{bound-3} \| \Phi - A_1 \Phi_1 - A_2 \Phi_2 - A_3 \Phi_3
\|_{C^0_b(\mathbb{R}^2)} \leq C \epsilon^{\tilde{r}}, \qquad
\tilde{r} = \min(4r-1,1-2r),
\end{equation}
where $\frac{1}{4} < r < \frac{1}{2}$. The fastest convergence
rate is $\epsilon^{1/3}$, which is obtained at $r=\frac{1}{3}$.
Consequently, the solution $\phi(x) = \sqrt{\epsilon}\Phi(x)$ of
the nonlinear elliptic problem \eqref{stationary} is
$\mathcal{O}(\epsilon^{\tilde{r}+1/2})$ accurate to the
approximation $\phi^{(1)}(x) = \sqrt{\epsilon}(A_1 \Phi_1 + A_2
\Phi_2 + A_3 \Phi_3)$, with the fastest convergence rate being $\epsilon^{5/6}$.
\end{corollary}

\begin{proof}
Using the triangle inequality and the bound \eqref{bound-2}, we
obtain
$$
\| \Phi - A_1 \Phi_1 - A_2 \Phi_2 - A_3 \Phi_3
\|_{C^0_b(\mathbb{R}^2)} \leq C_1\epsilon^{1-2r} +
C_2\sum_{j=1}^3\|A_j-B_j\|_{C^0_b(\mathbb{R}^2)},
$$
for some $C_1,C_2 > 0$. By the H\"{o}lder inequality we have
$$
\|A_j-B_j\|_{C^0_b(\mathbb{R}^2)}\leq
\|\hat{A}_j-\hat{B}_j\|_{L^1(\mathbb{R}^2)} =
\|\hat{A}_j-\hat{B}_j\|_{L^1(D_\epsilon)} + \| \hat{A}_j
\|_{L^1(\mathbb{R}^2\backslash D_\epsilon)},
$$
with the last equality due to a compact support of $\hat{B}_{1,2,3}(p)$ on $p \in
D_{\epsilon}$. The first term in the upper bound is
$\mathcal{O}(\epsilon^{\tilde{r}})$ by the bound
\eqref{bound-final} and the second term is smaller than any power
of $\epsilon$ if $\hat{\bf A} \in
L^1_q(\mathbb{R}^2,\mathbb{C}^3)$ for all $q \geq 0$.
\end{proof}

\section{Numerical approximations of localized reversible solutions}
\label{S:numerical}

We approximate here several classes of localized reversible
solutions of the differential coupled-mode system
\eqref{E:CME_diff} numerically. We use the same potential as in
the Examples \ref{example-1}--\ref{example-3}, i.e. $W(x) =
1-\cos(x)$ so that the bifurcation takes place for $\eta_0 \approx
0.1745$. For selected representative cases we also verify the convergence rate of Corollary
\ref{C:conv_rate_final} and the
conditions of Theorem \ref{theorem-3} on the Jacobian operator,
which guarantees persistence of these solutions in the full system
\eqref{stationary}. We limit our attention to the following
classes of reversible solutions:

\begin{itemize}
\item[(A)] defocusing case $\sigma=1$: $A_1=A_2=0$,
$A_3=R(r)e^{im\theta}$, $r=\frac{1}{\sqrt{\alpha_3}}\sqrt{y_1^2+y_2^2}$, $\theta =
\text{arg}(y_1 + i y_2)$

$m=0$ \ldots radially symmetric positive soliton\\
$m=1$ \ldots vortex of charge one

\item[(B)] focusing case $\sigma =-1$: $\ A_3=0$
\begin{itemize}
\item[(i)] $A_2=0$, $A_1=R(r)e^{im\theta}$,
$r=\sqrt{\frac{y_1^2}{\alpha_1}+\frac{y_2^2}{\alpha_2}}$, $\theta
= \text{arg}(\alpha_2 y_1+i\alpha_1 y_2)$

$m=0$ \ldots ellipsoidal positive soliton\\
$m=1$ \ldots ellipsoidal vortex of charge one

\item[(ii)] $A_1(y_1,y_2)=\pm A_2(y_2,y_1)\in\mathbb{R}, \quad A_1(y_1,y_2)=A_1(-y_1,y_2)=A_1(y_1,-y_2)$ \\
\ldots symmetric real coupled soliton
\item[(iii)] $A_1(y_1,y_2)=\pm i A_2(y_2,y_1) \in i\mathbb{R},  \quad A_1(y_1,y_2)=A_1(-y_1,y_2)=A_1(y_1,-y_2)$\\
\hspace{5cm}  \ldots $\pi/2$-phase delay coupled soliton
\item[(iv)] $A_1(y_1,y_2) = \pm i \bar{A_2}(y_2,y_1) \in \mathbb{C}, \quad A_1(y_1,y_2)=-\bar{A_1}(-y_1,y_2)=\bar{A_1}(y_1,-y_2)$\\
\hspace{5cm}  \ldots coupled vortex of charge one
\end{itemize}
\end{itemize}

We will see that the localized solutions with $A_3=0$ bifurcate at $\sigma =
-1$ from the upper edge $\Omega=\beta_1$ of the band gap
and the localized solutions with
$A_1=A_2=0$ bifurcate at $\sigma = +1$ from the lower edge $\Omega=\beta_2$
of the band gap.
All solutions above satisfy the reversibility condition of
Definition \ref{definition-reversible}. Provided the assumptions
of Theorem \ref{theorem-3} are satisfied, Corollary
\ref{C:reversibility} then guarantees that these solutions
correspond to reversible solutions of the nonlinear elliptic
problem \eqref{stationary}. Using the symmetries of $\psi_1$,
$\varphi_1$ and $\varphi_2$, we can see that the function
$\phi(x)$ satisfies the following symmetry:
\begin{itemize}
\item[(A)]
\begin{tabular}{ll}
$\phi(x_1,x_2)=\phi(-x_1,x_2)=\phi(x_1,-x_2) \in \mathbb{R}$ & $(m=0)$\\
$\phi(x_1,x_2)=-\bar{\phi}(-x_1,x_2)=\bar{\phi}(x_1,-x_2) \in
\mathbb{C}$ & $(m=1)$
\end{tabular}
\item[(B-i)]
\begin{tabular}{ll}
$\phi(x_1,x_2)=\phi(-x_1,x_2)=-\phi(x_1,-x_2) \in \mathbb{R}$ & ($m=0$)\\
$\phi(x_1,x_2)=-\bar{\phi}(-x_1,x_2)=-\bar{\phi}(x_1,-x_2) \in
\mathbb{C}$ & ($m=1$)
\end{tabular}
\item[(B-ii)] \begin{tabular}{l} $\phi(x_1,x_2)=\pm
\phi(x_2,x_1)=\mp\phi(-x_2,-x_1) \in \mathbb{R}$\end{tabular}
\item[(B-iii)]
\begin{tabular}{l} $\phi(x_1,x_2)=\pm i \bar{\phi}(x_2,x_1)=\mp i
\bar{\phi}(-x_2,-x_1) \in \mathbb{C}$\end{tabular}
\item[(B-iv)]
\begin{tabular}{l} $\phi(x_1,x_2)=\pm i \bar{\phi}(x_2,x_1)=\pm i
\bar{\phi}(-x_2,-x_1) \in \mathbb{C}$\end{tabular}
\end{itemize}
Note that (B-iii) and (B-iv) agree with the reversibility
constraint (\ref{reversibility-2}) after multiplication by
$e^{-i\pi/4}$.

The above solutions do not include any three-component gap
solitons with all $A_1,A_2$ and $A_3$ nonzero. Such solutions do
not bifurcate from either upper or lower edge $\Omega=\beta_1$ or
$\Omega=\beta_2$, respectively, and it is thus hard to capture such solutions
numerically. It
is, nevertheless, possible that the three-component solitons still
exist in the interior of the gap $\beta_2 < \Omega < \beta_1$ but
we do not attempt here to locate these special solutions.

\subsection{One-component solutions}\label{S:1-comp}

One-component solutions correspond to classes (A) and (B-i). The
function $R(r)$ for class (A) satisfies the ODE
\begin{equation}
R''+\frac{1}{r}R'+(\Omega-\beta_2)R-\frac{m^2}{r^2}R-\sigma
\gamma_4 R^3=0, \label{ODE-radial}
\end{equation}
where $R(0) > 0$, $R'(0)=0$ for $m=0$ and $R(0) = 0$, $R'(0) > 0$
for $m=1$. For $m\neq 0$ the initial-value problem for the ODE
(\ref{ODE-radial}) is ill-posed but can be turned into a
well-posed one via the transformation $Q=r^{-m} R(r)$ leading to
\begin{equation}\label{E:Q_equ}
Q''+\frac{2m+1}{r}Q'+(\Omega-\beta_2)Q-\sigma \gamma_4 r^{2m}
Q^3=0,
\end{equation}
with $Q(0) > 0$, such that $R(r) \sim r^{|m|}$ as $r\rightarrow
0$. Similarly, the function $R(r)$ for class (B-i) satisfies the
ODE
\begin{equation}\label{E:Q_equ2}
R''+\frac{1}{r} R' + (\Omega-\beta_1) R - \frac{m^2}{r^2} R -
\sigma \gamma_1 R^3=0.
\end{equation}
We solve these equations numerically via a shooting method
searching for $R(r)$ vanishing at infinity as $r\rightarrow
\infty$.

Figure \ref{F:cont_lines_1_comp} shows  the solution families in
the frequency-amplitude space for both solutions with $m = 0$ and
$m = 1$. Class (A) families bifurcate from the lower edge $\Omega = \beta_2$ and class (B-i)
families from the upper edge $\Omega = \beta_1$. Examples of the
gap solitons $(m=0)$ are in Figure \ref{F:single_comp_townes} and
of the vortices ($m=1$) in Figure \ref{F:single_comp_vort}. As
expected, the class (A) solutions are radially symmetric and the
class (B) solutions ellipsoidal.

For the case of one-component solutions, the coupled-mode
equations reduce to a scalar two-dimensional NLS equation.
Conditions of Theorem \ref{theorem-3} on the kernel of the
Jacobian operator are known to be satisfied for $m = 0$
\cite{Weinstein} and, therefore, do not need to be checked
numerically.

Figure \ref{F:eps_conv_1_comp} presents the computed
$\epsilon$-convergence of the error term $\phi(x) - \phi^{(1)}(x)$
described by Corollary \ref{C:conv_rate_final} for the solution classes (A)
and (B-i) with $m=0$. The convergence rate is seen
around $\epsilon^{1.07}$ and $\epsilon^{1.08}$, which is higher
than the fastest convergence rate $\epsilon^{5/6}$ predicted by Corollary
\ref{C:conv_rate_final}. The observed rate is close to $\epsilon^{1}$,
which is the  rate predicted by a formal multiple scales asymptotic expansion
of $\phi(x)$, for which $\phi^{(1)}(x)$ is the leading order term.

\begin{figure}[h!]
\begin{center}
\includegraphics[height=6cm]{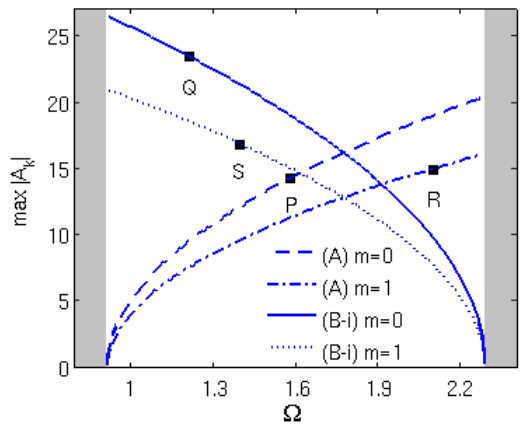}
\end{center}
\caption{Continuation curves of one-component solutions (A) and
(B-i). The marked points P,Q, R and S at $\Omega \approx 1.6,
1.22, 2.1$ and $1.4$ correspond to the profiles in Fig.
\ref{F:single_comp_townes} and \ref{F:single_comp_vort}.}
\label{F:cont_lines_1_comp}
\end{figure}

\begin{figure}[h!]
\begin{center}
\includegraphics[height=7.8cm]{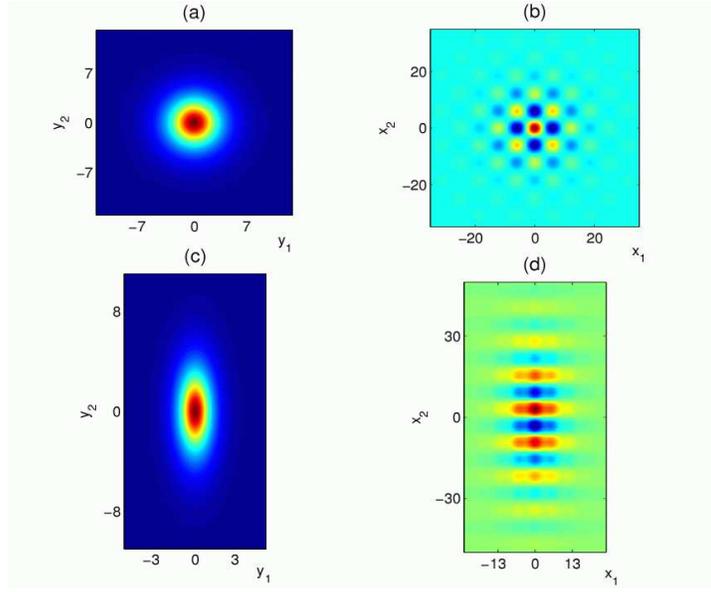}
\end{center}
\caption{Profiles of the one-component real gap soliton $(m = 0)$. (a)
$A_3$ at point P in Fig. \ref{F:cont_lines_1_comp};
(b) the corresponding leading-order term $A_3(y_1,y_2)\Phi_3(x_1,x_2)$ for
$\epsilon =0.1$; (c)  $A_1$ at point Q in Fig.
\ref{F:cont_lines_1_comp}; (d) the corresponding leading-order term
$A_1(y_1,y_2)\Phi_1(x_1,x_2)$ for $\epsilon=0.02$.}
\label{F:single_comp_townes}
\end{figure}

\begin{figure}[h!]
\begin{center}
\includegraphics[height=7.8cm]{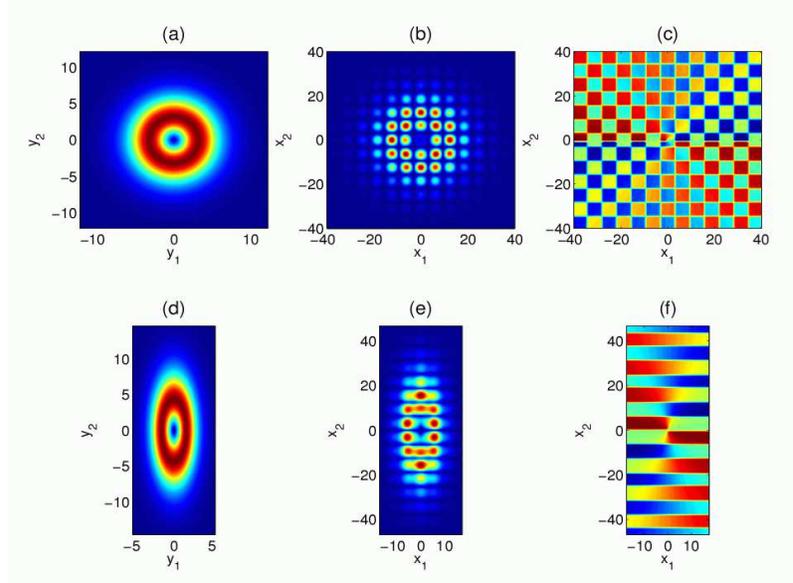}
\end{center}
\caption{Profiles of the one-component vortex ($m = 1$). (a)
$|A_3|$ at point R in Fig. \ref{F:cont_lines_1_comp};
(b) and (c) the modulus and phase, respectively, of the corresponding
leading-order term $A_3(y_1,y_2)\Phi_3(x_1,x_2)$ for $\epsilon =0.1$; (d)
$|A_1|$ at point S in Fig. \ref{F:cont_lines_1_comp};
(e) and (f) the modulus and phase, respectively, of the corresponding
leading-order term $A_1(y_1,y_2)\Phi_1(x_1,x_2)$ for $\epsilon=0.1$. }
\label{F:single_comp_vort}
\end{figure}

\begin{figure}[h!]
\begin{center}
\includegraphics[height=4cm]{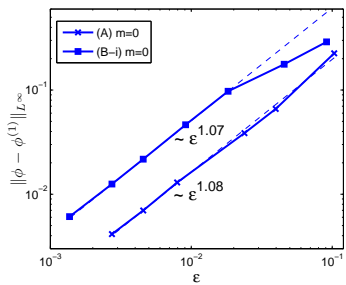}
\end{center}
\caption{$\epsilon$-convergence of the error term for the gap
soliton ($m = 0$) corresponding to (A) at $\Omega \approx 1.22$
and (B-i) at $\Omega \approx 1.25$.} \label{F:eps_conv_1_comp}
\end{figure}

\subsection{Two-component solutions}
\label{S:2-comp}

The following algorithm is used to compute two-component solutions
of classes (B-ii), (B-iii) and (B-iv). First, $\alpha_1$ and
$\alpha_2$ are replaced by
$\bar{\alpha}=\frac{1}{2}(\alpha_1+\alpha_2)$ so that the first
two equations of the coupled-mode system \eqref{E:CME_diff} admit
a reduction $A_1=A_2$, where $A_1=R(r)e^{im\theta}$ solves a
one-component problem. The function $R(r)$ is then computed as in
Section \ref{S:1-comp}. Next, homotopy in the coefficients
$\alpha_1$ and $\alpha_2$ is employed to gradually deform the
computed solution with $\alpha_1=\alpha_2=\bar{\alpha}$ into the
one with the original values of $\alpha_1$ and $\alpha_2$. In this
homotopy continuation the full coupled-mode system needs to be
solved numerically and this is done via Newton's method on a
central finite-difference discretization with zero Dirichlet
boundary conditions.

Figure \ref{F:cont_lines_2_comp} shows the solution families of
classes (B-ii), (B-iii) and (B-iv) as curves in the
frequency-amplitude space. Profiles of examples of these solutions
corresponding to the marked points in Figure
\ref{F:cont_lines_2_comp} appear in Figures \ref{F:two_comp_real},
\ref{F:two_comp_pi_half_delay} and \ref{F:two_comp_vort}. Note
that the coupled vortex in Figure \ref{F:two_comp_vort} has been
obtained via the above described homotopy continuation from a
vortex of charge one with $\alpha_1=\alpha_2$. As the phase plots
of $A_1$ and $A_2$ show, the resulting coupled vortex is also of
charge one.

In order to verify the persistence conditions of Theorem
\ref{theorem-3} for the two-component solutions, the kernel of the
Jacobian operator $J$ has to be studied numerically. When
rewritten in the real variables (for real and imaginary part) the
first two equations of the coupled-mode system \eqref{E:CME_diff}
give rise to a Jacobian operator $J$. For the cases (B-ii) and
(B-iii), the Jacobian operator can be block-diagonalized to the
form
\[
\begin{pmatrix} J_+ & 0\\ 0 & J_-\end{pmatrix}
\]
such that $\text{dim}(\text{Ker}(J)) =
\text{dim}(\text{Ker}(J_+))+\text{dim}(\text{Ker}(J_-))$. The
elements of $\text{Ker}(J)$ can, moreover, be easily reconstructed
from the elements of $\text{Ker}(J_+)$ and $\text{Ker}(J_-)$.
This diagonalization reduces the expense of the eigenvalue
solver.

\begin{figure}[h!]
\begin{center}
\includegraphics[height=3.5cm]{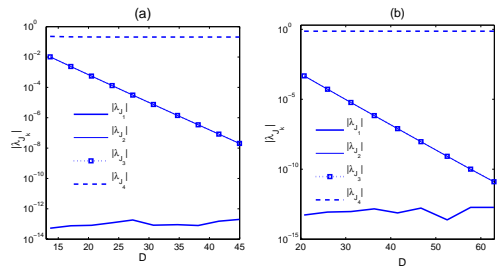}
\end{center}
\caption{The four smallest eigenvalues of the Jacobian operator
$J$ for (a) the two-component real gap soliton (B-ii) and (b) the
two-component $\pi/2$-phase delay gap soliton (B-iii) for $\Omega
\approx 1.19$.} \label{F:kernel_2_comp_real_delay}
\end{figure}

Figure \ref{F:kernel_2_comp_real_delay} shows the four smallest
(in modulus) eigenvalues $\lambda_{J_1}, \ldots, \lambda_{J_4}$ of
the 4th order central finite-difference approximation of $J$ for
examples of solutions (B-ii) and (B-iii) as functions of $D$, where
$D$ is the size
of the computational box $[-D,D]\times [-D,D] \subset \mathbb{R}^2$
with the step sizes $dy_1=dy_2\approx 0.14$ and with zero Dirichlet boundary conditions. Clearly, the three
smallest eigenvalues converge to zero as $D$ grows while the
fourth eigenvalue converges to a nonzero value. Therefore, the
assumptions of Theorem \ref{theorem-3} are verified for these
solutions. We have also checked that the corresponding
eigenvectors of the three zero eigenvalues are approximations of
$\partial_{y_1} {\bf A},\partial_{y_2} {\bf A}$ and $i {\bf A}$.

\begin{figure}[h!]
\begin{center}
\includegraphics[height=3.75cm]{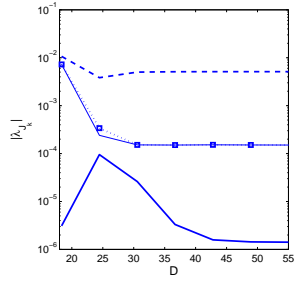}
\end{center}
\caption{The four smallest eigenvalues of the Jacobian operator
$J$ for the two-component vortex (B-iv) at $\Omega \approx 1.2$.}
\label{F:kernel_2_comp_vort}
\end{figure}

For the case (B-iv) the Jacobian operator $J$ has to be treated in
full because it cannot be block-diagonalized for complex-valued
solutions. Figure \ref{F:kernel_2_comp_vort} shows the four
smallest eigenvalues of the discretized Jacobian operator for an
example of the solution (B-iv). The fourth-order central
finite-difference discretization was used once again with the step
sizes $dy_1=dy_2 \approx 0.12$. The fourth eigenvalue
$\lambda_{J_4}$ converges approximately to $0.005$ as $D$ grows,
while the second and third eigenvalues $\lambda_{J_{2,3}}$
converge to $0.0002$.
We have checked that the limit of $\lambda_{J_{2,3}}$ decreases as
$dy_1=dy_2$ is decreased while the limit of $\lambda_{J_4}$
remains practically unchanged.

\begin{figure}[h!]
\begin{center}
\includegraphics[height=4cm]{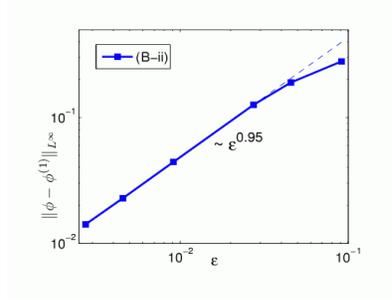}
\end{center}
\caption{$\epsilon$-convergence of the error term for the gap
soliton (B-ii) with $\Omega \approx 0.94$.}
\label{F:eps_conv_2_comp}
\end{figure}

Finally, Figure \ref{F:eps_conv_2_comp} verifies the
$\epsilon$-convergence of the error term $\phi(x) - \phi^{(1)}(x)$
described by Corollary \ref{C:conv_rate_final} for the gap soliton
of class (B-ii) at $\Omega \approx 0.944$. The observed
convergence rate is $\epsilon^{0.95}$, which is, once again, close
to $\epsilon^{1}$ expected from the formal asymptotic expansion.
The $\epsilon$-convergence has not
been checked for any complex-valued solutions of classes (B-iii)
and (B-iv) due to four times larger memory requirements compared
to the real case. The memory use grows rapidly as $\epsilon$
decreases since domains of the size $(2D/\sqrt{\epsilon}) \times
(2D/\sqrt{\epsilon})$ need to be discretized to compute the
solution $\phi(x)$ of the elliptic problem \eqref{stationary} and
each period $[0,2\pi] \times [0,2\pi]$ of the potential function
$V(x)$ requires about 100 grid points.

\begin{figure}[htbp]
\begin{center}
\includegraphics[height=4cm]{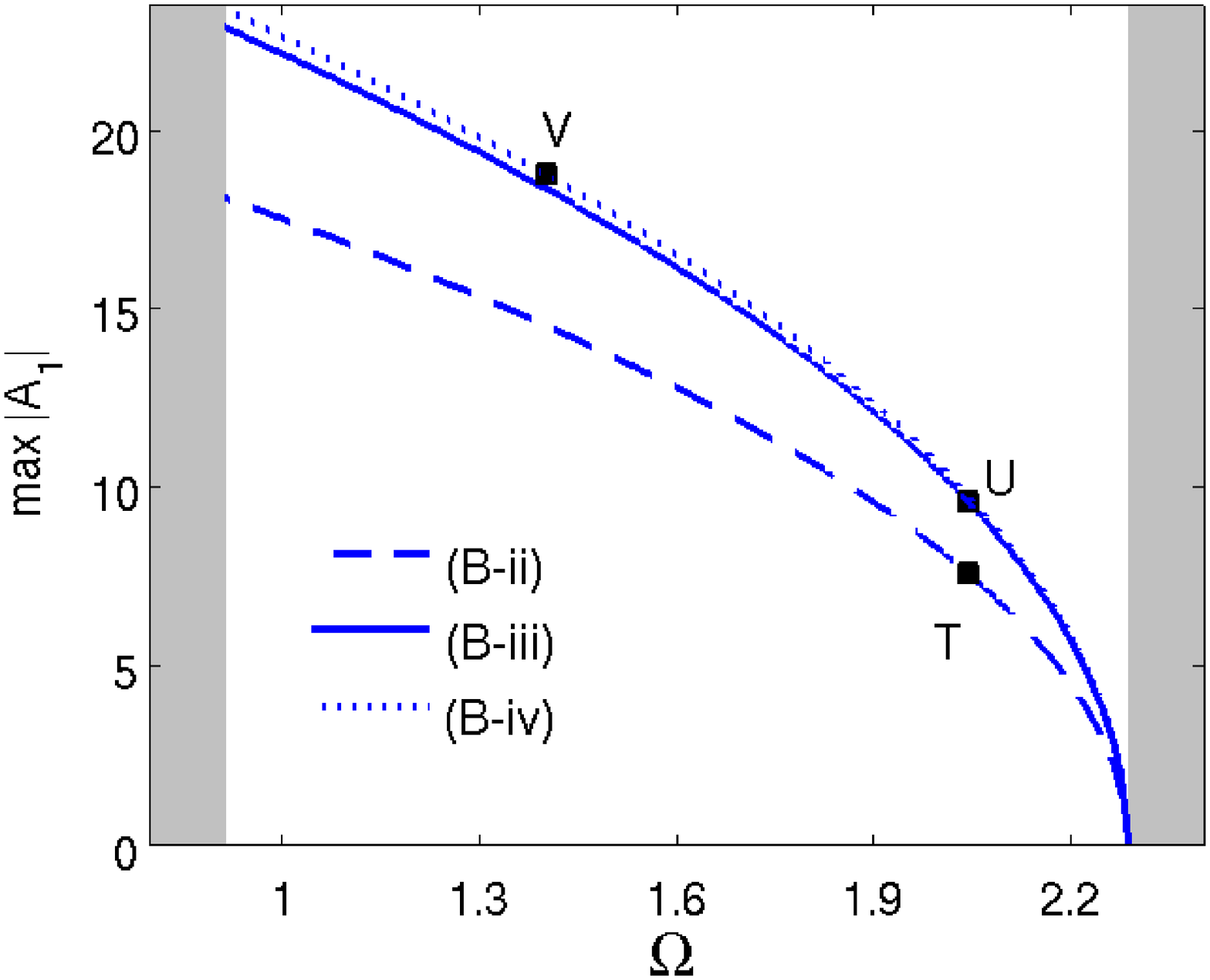}
\end{center}
\caption{Continuation curves of two-component solutions (B-ii),
(B-iii) and (B-iv). The marked points T and U at $\Omega \approx
2.04$ correspond to the profiles in Figs. \ref{F:two_comp_real}
and \ref{F:two_comp_pi_half_delay}, while the point V at $\Omega
\approx 1.4$ corresponds to the profiles in Fig.
\ref{F:two_comp_vort}.} \label{F:cont_lines_2_comp}
\end{figure}

\begin{figure}[h!]
\begin{center}
\includegraphics[height=4.75cm]{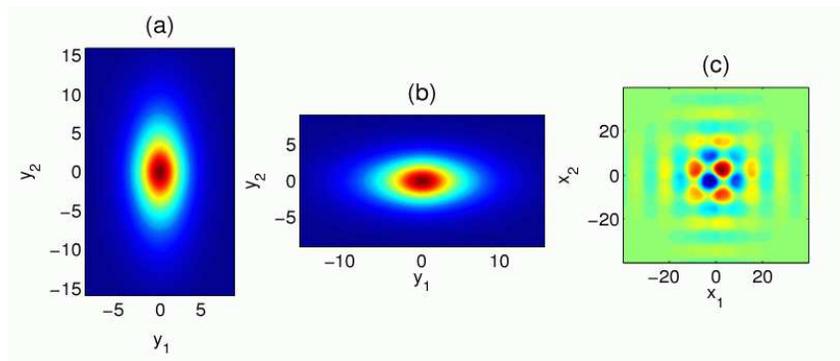}
\end{center}
\caption{Two-component real gap soliton (B-ii) at point T
in Fig. \ref{F:cont_lines_2_comp}. (a) $A_1$; (b) $A_2$; (c) the
leading-order term
$A_1(y_1,y_2)\Phi_1(x_1,x_2)+A_2(y_1,y_2)\Phi_2(x_1,x_2)$ for
$\epsilon=0.1$.} \label{F:two_comp_real}
\end{figure}

\begin{figure}[h!]
\begin{center}
\includegraphics[height=4cm]{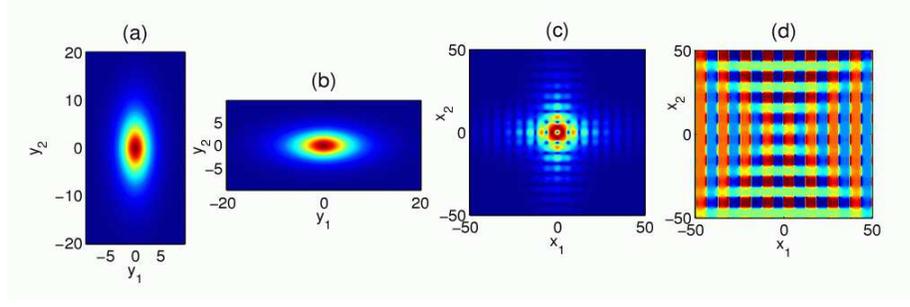}
\end{center}
\caption{Two-component $\pi/2$-phase delay gap soliton (B-iii)
at point U in Fig. \ref{F:cont_lines_2_comp}. (a)
$A_1$; (b) $-iA_2$; (c) and (d)  modulus and phase, respectively,
of the leading-order term
$A_1(y_1,y_2)\Phi_1(x_1,x_2)+A_2(y_1,y_2)\Phi_2(x_1,x_2)$ for
$\epsilon=0.1$.} \label{F:two_comp_pi_half_delay}.
\end{figure}

\begin{figure}[h!]
\begin{center}
\includegraphics[height=6cm]{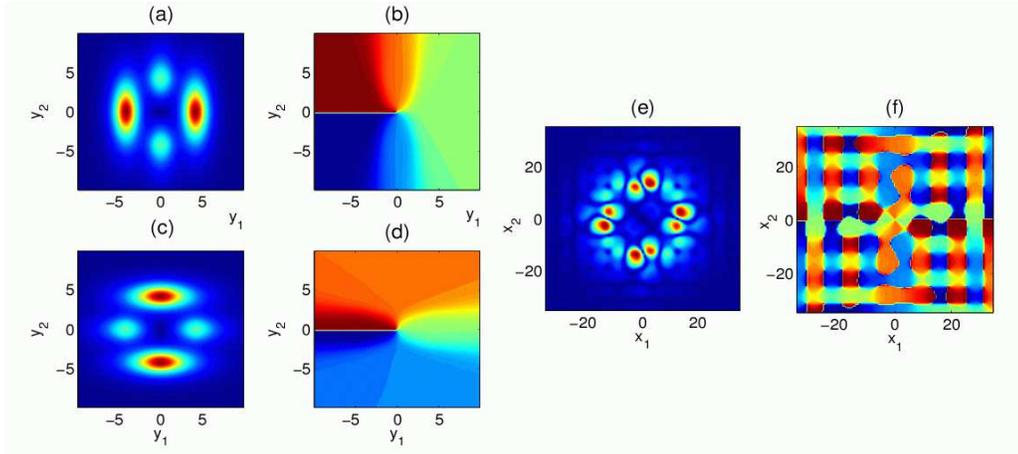}
\end{center}
\caption{Two-component vortex (B-iv) at point V in
Fig. \ref{F:cont_lines_2_comp}. (a) and (b) modulus and phase of
$A_1$; (c) and (d) modulus and phase of $A_2$; (e) and (f) modulus
and phase of the leading-order term
$A_1(y_1,y_2)\Phi_1(x_1,x_2)+A_2(y_1,y_2)\Phi_2(x_1,x_2)$ for
$\epsilon=0.1$.} \label{F:two_comp_vort}
\end{figure}

\section{The time-dependent case}\label{S:t_dep_CME}

We show here that dynamics of non-stationary localized solutions
of the Gross--Pitaevskii equation (\ref{GP}) are close to the
dynamics of the time-dependent coupled-mode equations for finite
time intervals. These results are similar to the ones in
\cite{BSTU,SU} and they give a rigorous basis for formal
asymptotic results in \cite{Konotop,SY}. According to a formal
asymptotic multiple scales expansion, solutions of the Gross--Pitaevskii equation \eqref{GP}
in the form
\begin{eqnarray*}
E_{\rm ans} (x,t) = \sqrt{\epsilon} E_1(x,t) +
\mathcal{O}(\epsilon),
\end{eqnarray*}
where
\begin{eqnarray*}
E_1 = \left[ A_1 (\sqrt{\epsilon}x, \epsilon t) \psi_1 (x_1)
\varphi_2(x_2) + A_2 (\sqrt{\epsilon}x, \epsilon t) \varphi_2
(x_1) \psi_1(x_2) + A_3 (\sqrt{\epsilon} x, \epsilon t) \varphi_1
(x_1) \varphi_1 (x_2) \right] e^{i\omega_0 t},
\end{eqnarray*}
are approximated by solutions of the time-dependent coupled-mode
equations
{\small
\begin{equation}
\left\{
\begin{array}{rcl} (i\partial_T - \beta_1) A_1 + \left( \alpha_1
\partial_{y_1}^2 + \alpha_2 \partial_{y_2}^2 \right) A_1 & = &
\sigma \left[ \gamma_1 |A_1|^2 A_1 + \gamma_2 (2 |A_2|^2 A_1 +
A_2^2 \bar{A}_1)
+ \gamma_3 (2 |A_3|^2 A_1 + A_3^2 \bar{A}_1) \right], \\
(i\partial_T - \beta_1) A_2 + \left( \alpha_2 \partial_{y_1}^2 +
\alpha_1 \partial_{y_2}^2 \right) A_2 & = & \sigma \left[ \gamma_1
|A_2|^2 A_2 + \gamma_2 (2 |A_1|^2 A_2 + A_1^2 \bar{A}_2)
+ \gamma_3 (2 |A_3|^2 A_2 + A_3^2 \bar{A}_2) \right], \\
(i \partial_T - \beta_2) A_3 + \alpha_3 \left( \partial_{y_1}^2 +
\partial_{y_2}^2 \right) A_3 & = & \sigma \left[ \gamma_4 |A_3|^2
A_3 + 2 \gamma_3 (|A_1|^2 + |A_2|^2) A_3 + \gamma_3 (A_1^2 +
A_2^2) \bar{A}_3 \right],
\end{array}  \right. \label{time-dependent-cme}
\end{equation}
}where $y = \sqrt{\epsilon} x$ and $T = \epsilon t$. Under the
non-resonance condition
\begin{equation} \label{gs4}
\inf_{{(n_1,n_2) \in \mathbb{N}^2} \atop {j_1, j_2, j_3 \in \{\pm
5, \pm 3, \pm 1\}}  } \left| \rho_{n_1}\left(\frac{j_2 +
j_3}{2}\right)+ \rho_{n_2}\left(\frac{j_1 + j_3}{2}\right) - |j_1 +
j_2 + j_3| \omega_0 \right| > 0,
\end{equation}
where $|j_1 + j_2 + j_3 | \leq 5$ and the infinimum does not
include direct resonances that give the nonlinear terms in system (\ref{time-dependent-cme}),
the following theorem states the approximation result on the time-dependent solutions.

\begin{theorem}
Let $s > 1, s_* \geq \max\{3,s\}$, and assume that the non-resonance
condition \eqref{gs4} holds. Then for all $C_1$ and $T_0
> 0$ there exist $\epsilon_0 > 0$ and $C_2 > 0$ such that for
all solutions $A_1, A_2, A_3 \in C([0, T_0], H^{s_*}(\mathbb{R}^2))$ of the
time-dependent coupled-mode system with
\[
\sup_{T\in [0, T_0]} \|{\bf A}(\cdot, T)\|_{H^{s_*}(\mathbb{R}^2)}
\leq C_1,
\]
there exist solutions $E\in C([0, T_0/\epsilon], H^s(\mathbb{R}^2))$ of the
Gross--Pitaevskii equation \eqref{GP} with
\[
\sup_{t\in [0, T_0/\epsilon] } \|E (\cdot, t) - E_{ans}(\cdot,
t)\|_{H^s(\mathbb{R}^2)} \leq C_2\epsilon^{3/2}.
\]
for all $\epsilon \in (0, \epsilon_0)$.
\end{theorem}

\begin{proof}
The proof is very similar to the one given for a scalar nonlinear
Schr\"{o}dinger equation in the one-dimensional case in
\cite{BSTU}. The non-resonance condition (\ref{gs4}) differs from
\cite[Eq.(4)]{BSTU} by including the terms with $j_1,j_2,j_3 = \pm
5$. This is necessary to make the residual formally of order
$\mathcal{O}(\epsilon^3)$. Due to the fact that we are working
here on $x \in \mathbb{R}^2$, we lose $\epsilon^{-1/2}$ due to
the scaling of the $L^2$-norm in contrast to \cite{BSTU}, where
we would only lose $\epsilon^{-1/4}$ on $x \in \mathbb{R}$ with
the same scaling. Assuming that the non-resonance condition
\eqref{gs4} is satisfied, the $\mathcal{O}(\epsilon)$ terms in
$E_{\rm ans}$ can be chosen in such a way that all terms up to order
$\mathcal{O} (\epsilon^{5/2})$ in the residual
\[
\mathop{\mathrm{Res}} (E_{\rm ans}) = - i E_t - \nabla^2 E + V(x) E +
\sigma |E|^2 E
\]
can be eliminated. By \cite[Lemma 3.3]{BSTU} Bloch transform is an
isomorphism between $H^s(\mathbb{R}^2)$ and $l^2_s (\mathbb{N}^2,
L^2(\mathbb{T}^2))$. Hence, similarly to \cite[Lemma 3.4]{BSTU} the
nonlinearity $E \mapsto |E|^2 E$ maps $l^2_s (\mathbb{N}^2,
L^2(\mathbb{T}^2))$ in Bloch space, or $H^s(\mathbb{R}^2)$ in
physical space, into itself for $s > 1$ due to Sobolev's embedding
theorem in two dimensions. Moreover, $L = - \nabla^2 + V(x)$
generates a uniformly bounded semigroup in $l^2_s (\mathbb{N}^2,
L^2(\mathbb{T}^2))$, or respectively in $H^s(\mathbb{R}^2)$,
according to \cite[page 927]{BSTU}. Hence, the error function $R =
\epsilon^{-3/2} (E-E_{ans})$ satisfies an equation of the form
\[
\partial_tR = L R + \mathcal{O} (\epsilon).
\]
Using the variation of constant formula and Gronwall's inequality
gives the required $\mathcal{O}(1)$ bound for $R$. For more details
see \cite[Section 4.2]{BSTU}.
\end{proof}

\section{Various generalizations}\label{S:generalizations}

We have given a detailed justification of the three
coupled-mode equations which describe gap solitons bifurcating in
the first band gap of a two-dimensional separable periodic
potential. These equations are generic (structurally stable) for
the considered bifurcation but they can be modified for other
relevant bifurcation problems. We review here several examples of
other bifurcations, where the justification analysis is expected
to be applicable in a similar manner.

{\em Bifurcations from band edges.} If the band edge separates the
existing band gap of finite length and a single or
double-degenerate spectral band, the coupled-mode equations
persist without any modifications. A formal derivation is given in
the recent paper \cite{SY}, where a scalar NLS equation is derived for
a single band and two coupled NLS equations are derived for a
double band. Note that the coupled NLS equations for a double band
may have terms destroying the reduction $A_1 = 0$ or $A_2 = 0$ in
the truncated system of equations. This happens ,for instance,at the
band edge $E$ in the second band gap of a separable symmetric
potential \cite{SY}. Since the analysis of separable potentials relies
on the construction of one-dimensional Bloch modes and
one-dimensional spectral bands, the same analysis is valid for
bifurcations from band edges in one-dimensional problems. The
formal derivation of the NLS equation for a band edge of a single
spectral band in the one-dimensional GP equation was reported in \cite{PSK}.

{\em Bifurcations in the higher-order band gaps.} With larger
values of the parameter $\eta$ in the separable potential
(\ref{separable-potential}), more band gaps open, although the
number of band gaps is always finite for finite $\eta$, see Theorem 6.10.5 in \cite{Eastham}. Since
the bifurcation of new band gaps occurs due to the resonance of
finitely many Bloch modes in the separable potential, a system of
finitely-many coupled-mode equations can be derived similarly
in the higher-order band gaps. For instance, the bifurcation
of the second band gap in the symmetric separable potential may
occur for a resonance of either three modes when $\lambda_2 +
\mu_1 = 2 \mu_2 < \lambda_1 + \lambda_3$ or four modes when
$\lambda_2 + \mu_1 = \lambda_1 + \lambda_3 > 2 \mu_2$ or five
modes when $\lambda_2 + \mu_1 = \lambda_1 + \lambda_3 = 2 \mu_2$.

{\em Bifurcations in asymmetric separable potentials.} If the
separable potential is asymmetric, such that
$$
V(x) = \eta \left( W_1(x_1) + W_2(x_2) \right),
$$
where $W_1$ and $W_2$ are two different functions on $x \in
\mathbb{R}$, the degeneracy of spectral bands is broken and
bifurcations of the band gaps occur with a smaller number of
resonant Bloch modes. Our analysis remains valid but the spectrum
of two different one-dimensional Sturm--Liouville problems must be
incorporated in the Fourier--Bloch decomposition. In this case,
two coupled-mode equations occur generally for the bifurcation of the
first band gap in the anisotropic separable potential.

{\em Bifurcations in finite-gap potentials.} Finite-gap potentials
lead to the degeneracy of eigenvalues of the one-dimensional
Sturm--Liouville problem and non-zero derivatives of the functions
$\rho_n(k)$ at the extremal points $k = 0$ and $k = \pm
\frac{1}{2}$. Adding a generic potential breaks the symmetry of
the finite-gap potential and leads to the bifurcation of narrow
band gaps in the space of one dimension. The corresponding
coupled-mode equations must have first-order than second-order derivative
operators. A
formal derivation of such coupled-mode equations from the Bloch
mode decomposition was considered in \cite{dSSS96}. Our
analysis provides a rigorous justification of these coupled-mode
equations.

{\em Bifurcation in super-lattices.} Let the potential $V(x)$ be
represented in the form
$$
V(x) = \eta_0 \left( W(x_1) + W(x_2) \right) + \epsilon
\left( \tilde{W}(x_1) + \tilde{W}(x_2) \right),
$$
where $W(x)$ and $\eta_0$ are the same as in the separable potential
(\ref{separable-potential}) and $\tilde{W}(x)$ is a $4 \pi$-periodic
potential. The perturbation term $\tilde{W}(x)$ couples resonant Bloch
modes of the potential $W(x)$ by linear terms and destroys
reductions $A_1 = 0$, $A_2 = 0$ and $A_3 = 0$ in the differential
coupled-mode system. The system is still formulated in the form of
the coupled NLS equations with various (linear and nonlinear) coupling terms.

{\em Bifurcations in three-dimensional periodic problems.} The
same analysis holds in three-dimensional separable potentials,
since the spectral bands are still enumerated by a countable
number of spectral bands of the one-dimensional potentials.
Bifurcations of new band gaps occur again due to a resonance of
finitely many Bloch modes of the one-dimensional potentials.

{\em Bifurcations in non-separable periodic potentials.} At the
present time, the generalization of the analysis for non-separable
periodic potentials meets a technical obstacle that the bound
(\ref{bound-K}) obtained in \cite{BSTU} is needed to be extended to
problems with two-dimensional periodic potentials.

{\bf Acknowledgement.} The second author thanks A. Sukhorukov for
discussions at the early stage of the project. The work of T. Dohnal and D.
Pelinovsky is supported by the Humboldt Research Fellowship. The
work of G. Schneider is partially supported  by the
Graduiertenkolleg 1294 ``Analysis, simulation and design of
nano-technological processes'' granted by the Deutsche
Forschungsgemeinschaft (DFG) and the Land Baden-W\"{u}rttemberg.


\begin{thebibliography}{99}

\bibitem{ACD95} A.B. Aceves, B. Costantini and C. De Angelis,
``Two-dimensional gap solitons in a nonlinear periodic slab
waveguide'',  J. Opt. Soc. Am. B {\bf 12}, 1475--1479 (1995)

\bibitem{AFI04} A.B. Aceves, G. Fibich and B. Ilan,
``Gap-soliton bullets in waveguide gratings'', Physica D {\bf
189}, 277--286 (2004)

\bibitem{AP05} D. Agueev and D. Pelinovsky, ``Modeling of wave
resonances in low-contrast photonic crystals'', SIAM J. Appl.
Math. {\bf 65}, 1101--1129 (2005)

\bibitem{AJ98} N. Ak\"{o}zbek and S. John,
``Optical solitary waves in two and three dimensional photonic
bandgap structures'', Phys. Rev. E {\bf 57}, 2287--2319 (1998)

\bibitem{arnold} V. Arnold, ``Remarks on perturbation theory for problems of Mathieu
type'', Uspekhi Mat. Nauk {\bf 38}, 189--203 (1983) [English translation:
Russian Math. Surveys {\bf 38}, 215--233 (1983)]

\bibitem{Konotop} V.A. Brazhnyi, V.V. Konotop, V. Kuzmiak, and V.S. Shchesnovich,
``Nonlinear tunneling in two-dimensional lattices'', Phys. Rev. A {\bf 76},
023608 (2007)

\bibitem{BSTU} K. Busch, G. Schneider, L. Tkeshelashvili, and H.
Uecker, ``Justification of the nonlinear Schr\"{o}dinger equation
in spatially periodic media'', Z. Angew. Math. Phys. {\bf 57},
905--939 (2006)

\bibitem{DA05} T. Dohnal and A.B. Aceves,
``Optical soliton bullets in (2+1)D nonlinear Bragg resonant
periodic geometries'', Stud. Appl. Math. {\bf 115}, 209--232
(2005)

\bibitem{dSSS96} C.M. de Sterke, D.G. Salinas and J.E. Sipe,
``Coupled-mode theory for light propagation through deep nonlinear
gratings'', Phys. Rev. E {\bf 54}, 1969--1989 (1996)

\bibitem{Eastham} M.S. Eastham, {\em The Spectral Theory of Periodic Differential Equations},
Scottish Academic Press, Edinburgh, 1973

\bibitem{gelfand} I.M. Gelfand, ``Expansion in eigenfunctions of an equation with periodic
coefficients'', Dokl. Akad. Nauk. SSSR {\bf 73} (1950), 1117-1120

\bibitem{GS84} M. Golubitsky and D.G. Schaeffer, {\em Singularities and
Groups in Bifurcation Theory}, v. 1 (Springer-Verlag, Berlin, 1985)

\bibitem{goodman} R.H. Goodman, M.I. Weinstein, and P.J. Holmes, ``Nonlinear propagation
of light in one-dimensional periodic structures'', J. Nonlinear.
Science {\bf 11}, 123--168 (2001)

\bibitem{HKSW_prep}
J.M. Harrison, P. Kuchment, A. Sobolev, and B. Winn, ``On
occurrence of spectral edges for periodic operators inside the
Brillouin zone'', J. Phys. A: Math. Theor. {\bf 40}, 7597-7618 (2007)

\bibitem{heinz} H.P. Heinz, T. K\"{u}pper, and C.A. Stuart, ``Existence
and bifurcation of solutions for nonlinear perturbations of the
periodic Schr\"{o}dinger equation'', J. Diff. Eqs. {\bf 100},
341--354 (1992)

\bibitem{Kuchment} P. Kuchment, ``The mathematics of photonic
crystals'', in {\em Mathematical Modeling in Optical Science}
(SIAM, Philadelphia, 2001)

\bibitem{kupper} T. K\"{u}pper and C.A. Stuart, ``Necessary and
sufficient conditions for gap-bifurcation'', Nonlin. Anal. {\bf
18}, 893--903 (1992)

\bibitem{mills} D.L. Mills, {\em Nonlinear Optics: Basic Concepts}
(Springer Verlag, 1984)

\bibitem{Pankov} A. Pankov, ``Periodic nonlinear Schr\"{o}dinger equation with
application to photonic crystals'', Milan J. Math. {\bf 73},
259--287 (2005)

\bibitem{Parnov_08} L. Parnovski, ``Bethe-Sommerfeld Conjecture,'' arXiv:0801.3096, 2008.

\bibitem{PSn07} D. Pelinovsky and G. Schneider, ``Justification of the coupled-mode
approximation for a nonlinear elliptic problem with a periodic
potential'', Appl. Anal. {\bf 86}, 1017--1036 (2007)

\bibitem{PSK} D.E. Pelinovsky, A.A. Sukhorukov, and Yu.S. Kivshar,
``Bifurcations and stability of gap solitons in periodic
potentials'', Phys. Rev. E {\bf 70}, 036618 (2004)

\bibitem{RS} M. Reed and B. Simon, {\em Methods of Modern
Mathematical Physics. IV. Analysis of Operators}, Academic Press,
New York, 1978

\bibitem{SU} G. Schneider and H. Uecker, ``Nonlinear coupled mode dynamics
in hyperbolic and parabolic periodically structured spatially
extended systems'', Asymp. Anal. {\bf 28}, 163--180 (2001)

\bibitem{SY} Z. Shi and J. Yang, ``Solitary waves bifurcated
from Bloch-band edges in two-dimensional periodic media'', Phys.
Rev. E {\bf 75}, 056602 (2007)

\bibitem{SS} C.M. de Sterke and J.E. Sipe, ``Gap solitons'',
Progr. Opt., {\bf 33}, 203 (1994)

\bibitem{Stuart} C.A. Stuart, ``Bifurcations into spectral gaps'', Bull. Belg. Math. Soc.
Simon Stevin, 1995, suppl., 59pp.

\bibitem{Weinstein} M. Weinstein, ``Lyapunov stability of ground states
of nonlinear dispersive evolution equations'', Comm. Pure Appl.
Math. {\bf 39}, 51--67 (1986)

\end{thebibliography}
\end{document}